\documentclass[11pt,letterpaper,twoside]{article}
\usepackage[T1]{fontenc}
\usepackage[utf8]{inputenc}
\usepackage{lmodern}
\usepackage[margin=1in,headheight=14pt]{geometry}
\usepackage{amsmath,amssymb,amsthm,mathtools}
\usepackage{microtype}
\usepackage{enumitem}
\usepackage{aliascnt}
\newif\ifshowreviewchanges
\showreviewchangestrue
\ifshowreviewchanges
  \usepackage{changes}
\else
  \usepackage[final]{changes}
\fi
\definechangesauthor[name={Author suggestion}, color=blue]{author}
\usepackage{fancyhdr}
\usepackage[hidelinks,bookmarksnumbered]{hyperref}
\usepackage[nameinlink,noabbrev]{cleveref}
\hypersetup{
  pdftitle={A direct obstacle proof of a dimension-free weak-type (1,1) bound for the Gaussian Riesz transforms},
  pdfauthor={Cristian E. Guti\'errez and Juan Manfredi},
  pdfsubject={Gaussian Riesz transforms, an exact obstacle decomposition, and dimension-free endpoint estimates},
  pdfkeywords={Gaussian measure, Ornstein-Uhlenbeck, Riesz transform, obstacle problem, weak type (1,1)}
}
\numberwithin{equation}{section}
\allowdisplaybreaks[2]
\newtheorem{theorem}{Theorem}[section]
\newaliascnt{proposition}{theorem}
\newtheorem{proposition}[proposition]{Proposition}
\aliascntresetthe{proposition}
\newaliascnt{lemma}{theorem}
\newtheorem{lemma}[lemma]{Lemma}
\aliascntresetthe{lemma}
\newaliascnt{corollary}{theorem}
\newtheorem{corollary}[corollary]{Corollary}
\aliascntresetthe{corollary}
\theoremstyle{definition}
\newaliascnt{definition}{theorem}
\newtheorem{definition}[definition]{Definition}
\aliascntresetthe{definition}
\theoremstyle{remark}
\newaliascnt{remark}{theorem}
\newtheorem{remark}[remark]{Remark}
\aliascntresetthe{remark}
\crefname{theorem}{Theorem}{Theorems}
\crefname{proposition}{Proposition}{Propositions}
\crefname{lemma}{Lemma}{Lemmas}
\crefname{corollary}{Corollary}{Corollaries}
\crefname{definition}{Definition}{Definitions}
\crefname{remark}{Remark}{Remarks}
\crefname{section}{Section}{Sections}
\crefname{subsection}{Section}{Sections}
\newcommand{\R}{\mathbb{R}}
\newcommand{\N}{\mathbb{N}}
\newcommand{\ga}{\gamma_n}
\newcommand{\dga}{\,d\gamma_n}
\newcommand{\Rg}{R_{\gamma}}
\newcommand{\E}{\mathcal{E}}
\newcommand{\U}{\mathsf{U}}
\newcommand{\V}{\mathsf{V}}
\newcommand{\K}{\mathsf{K}}
\newcommand{\Dom}{\operatorname{Dom}}

\newcommand{\one}{\mathbf{1}}
\newcommand{\ip}[2]{\left\langle #1,#2\right\rangle_2}
\newcommand{\dual}[2]{\left\langle #1,#2\right\rangle_{\V^*,\V}}
\newcommand{\norm}[1]{\left\lVert #1\right\rVert}
\newcommand{\abs}[1]{\left\lvert #1\right\rvert}

\newcommand{\me}{\mathfrak{m}}

\title{\vspace{-1.5em}\textbf{A direct obstacle proof of a dimension-free\\
weak-type $(1,1)$ bound for the Gaussian Riesz transforms}}
\author{Cristian E. Guti\'errez and Juan Manfredi}
\date{October 12, 2026}

\begin{document}
\maketitle
\thispagestyle{plain}

\begin{abstract} 
We give a detailed, direct proof of a dimension-free weak-type $(1,1)$
estimate for the first-order Gaussian vector Riesz transform, using the
fractional obstacle method of Ouyang, Spector, and Stockdale. For
$d\gamma_n=\pi^{-n/2}e^{-|x|^2}\,dx$, let
$L=-\tfrac12\Delta+x\cdot\nabla$, $D=2^{-1/2}\nabla$, and
$R_\gamma=DL^{-1/2}\Pi_0$, where $\Pi_0$ removes the mean. We prove that
$R_\gamma$ has weak-type constant at most $2$ on real-valued inputs,
independently of $n$. The main construction is the unique unshifted
decomposition $F=g+L^{1/2}u$, for nonnegative $F\in L^2(\gamma_n)$ and
$a>\int F\,d\gamma_n$, with $u\geq0$, $\min\{F,a\}\leq g\leq a$,
preservation of mass, and $g=a$ on $\{u>0\}$. We retain the details of
coercivity, the Lewy--Stampacchia estimate, the upgrade from the form
domain to the operator domain, and locality. We also identify the
$L^1$ extension with the principal-value transform and obtain the
constant $8/\pi$ for complex-valued inputs.
\end{abstract}



\tableofcontents

\section{Introduction and the main estimate}
\label{sec:context}

The first-order Riesz transforms associated with Gaussian measure are
of weak type $(1,1)$ in every finite dimension. In dimension one this
goes back to Muckenhoupt \cite{Muc69}. Fabes, Guti\'errez, and Scotto
\cite{FGS94} established the higher-dimensional theory, including weak-type
bounds for maximal spatial truncations, with constants depending on the
dimension. Further proofs of the fixed-dimensional endpoint estimates
were given by P\'erez and Soria \cite{PS00} and by Bruno \cite{Bruno19}.
The question of dimension-independent endpoint constants has a long
history: Guti\'errez \cite[p.~43]{Gut93} raised it in connection with the
maximal Gaussian transforms, while Ba\~nuelos formulated the question
for the untruncated Gaussian vector transform explicitly as
\cite[Problem~8]{Ban10}, explaining its connection with Wiener space.
The distinction between the linear and maximal questions is important;
our proof concerns the former.

In contrast, the dimension-free $L^p$ theory for $1<p<\infty$ has been
available since the work of Meyer \cite{Meyer84}; see Pisier
\cite{Pis88} and Guti\'errez \cite{Gut94} for analytic approaches.
These strong-type estimates do not supply the weak endpoint merely by
letting $p\downarrow1$. The argument below uses the exact $L^2$ energy
identity and an obstacle decomposition instead.

Ouyang, Spector, and Stockdale \cite{OSS26} obtain the weak-type
constant $2$ for the full Euclidean vector Riesz transform through a
fractional obstacle problem. Their remarkable decomposition is adapted
to the factorization of the operator, rather than to a covering by balls
or cubes. Mao, Wang, and Zhang \cite{MWZ26} extend this approach to
stratified Lie groups.

The Gaussian endpoint estimate is also covered by the abstract theorem
of Chen, Jiang, Li, and Li \cite[Theorem~5.1]{CJLL26}. In that theorem,
one takes $D=2^{-1/2}\nabla$ on $L^2(\gamma_n)$, so that $D^*D=L$;
the Ornstein--Uhlenbeck semigroup is Markovian and the Sobolev gradient
vanishes almost everywhere on a zero set. These hypotheses give the
same constant $2$. Their proof uses an obstacle problem for
$L^{1/2}+\kappa I$ and then removes the positive shift. Here we give a
self-contained proof directly for $L^{1/2}$, retaining exact mass
preservation and uniqueness of the obstacle decomposition.

Chen, Chewi, Lu, and Zhang \cite[Proposition~3.1 and Theorem~3.3]{CCLZ26}
use an unshifted obstacle problem for the Ornstein--Uhlenbeck operator
itself to obtain a second-order Gaussian endpoint estimate. Their
finite-mass argument likewise uses a height above the mean. Our
first-order construction instead uses the square root, whose form
domain is $\operatorname{Dom}(L^{1/4})$. The contribution of the present treatment is the direct, mass-preserving
square-root decomposition, with the intermediate functional-analytic
steps and normalizations made explicit.

We use the normalization
\[
 d\gamma_n=\pi^{-n/2}e^{-|x|^2}\,dx,\qquad
 L=-\tfrac12\Delta+x\cdot\nabla,\qquad
 D=2^{-1/2}\nabla,
\]
and set $R_\gamma=DL^{-1/2}\Pi_0$, with the inverse square root
annihilating constants. The structural identity is
\begin{equation}\label{eq:structural-intro}
 \Rg L^{1/2}u=Du,\qquad u\in\Dom(L^{1/2}).
\end{equation}
It converts the transformed remainder into the ordinary normalized
gradient of the obstacle potential. Arranging that this gradient
vanishes outside a set of controlled Gaussian measure reduces the
endpoint argument to an $L^2$ estimate for the good part.

The main endpoint estimate, proved in \cref{thm:endpoint}, is
\begin{equation}\label{eq:target}
 \sup_{t>0}t\,\ga\bigl(\{x:|\Rg f(x)|_{\R^n}>t\}\bigr)
 \leq2\norm{f}_{L^1(\ga)},\qquad f\in L^1(\ga;\R).
\end{equation}
Here $|\cdot|_{\R^n}$ is the Euclidean norm of the full vector, so the
same estimate holds for every coordinate. Throughout, for a measurable
vector-valued function $G$, we use the weak-$L^1$ quasi-norm
\[
 \norm{G}_{L^{1,\infty}(\ga)}
 :=\sup_{t>0}t\,\ga(\{|G|>t\}).
\]
The analogous notation for complex vectors uses their Euclidean norm.

The main construction is \cref{thm:decomposition}. For each
nonnegative $F\in L^2(\ga)$ and $a>\int F\dga$, it yields a unique
pair $u,g$ with $F=g+L^{1/2}u$, $u\geq0$, $\min\{F,a\}\leq g\leq a$,
and $\int g\dga=\int F\dga$. Moreover, $g=a$ on $\{u>0\}$ and
$Du$ vanishes on its complement. A decomposition with this bound on
$g$ cannot exist below the mean, since $L^{1/2}$ annihilates constants.
This is harmless for \eqref{eq:target}: on a probability space the
levels $t\leq\norm{f}_1$ are already trivial.

We keep the intermediate steps of the construction explicit.
\Cref{sec:spectral} fixes the normalization and the spectral domains.
\Cref{sec:decomposition} states the obstacle theorem and explains its
height restriction. In \cref{sec:proof}, the spectral gap controls the
mean-zero part of a minimizing sequence, while the difference between
the height and the mean controls its constant part. The first variation
gives a multiplier in the dual of $\Dom(L^{1/4})$. Its
Lewy--Stampacchia bound makes it an $L^2$ function and upgrades the
potential to $\Dom(L^{1/2})=W^{1,2}(\ga)$. Mass preservation,
complementarity, and the zero-set chain rule then finish the decomposition.
\Cref{sec:endpoint} proves the endpoint estimate, its extension to $L^1$,
and the principal-value identification. \Cref{app:pv-core} supplies the
kernel cancellation needed for that identification on smooth inputs.

\section{Gaussian normalization and spectral preliminaries}
\label{sec:spectral}

\subsection{The probability measure and the first-order energy}

Fix $n\geq1$ and set
\begin{equation}\label{eq:normalization}
  d\ga(x)=\pi^{-n/2}e^{-|x|^2}\,dx,
  \qquad L=-\frac12\Delta+x\cdot\nabla,
  \qquad D=\frac1{\sqrt2}\nabla.
\end{equation}
Unless complex-valued functions are explicitly considered, all spaces and
inner products are real and taken with respect to $\ga$; unmarked
integrals are over $\R^n$.  Since $\ga$ is a probability
measure, Cauchy--Schwarz gives
\begin{equation}\label{eq:L2intoL1}
  \norm{v}_1\leq\norm{v}_2.
\end{equation}

For $v,w\in C_c^\infty(\R^n)$, Gaussian integration by parts gives
\[
 \int (\partial_jv)w\dga
 =\int v(-\partial_jw+2x_jw)\dga.
\]
Hence
\[
 D_j^*=\frac1{\sqrt2}(-\partial_j+2x_j),
 \qquad \sum_{j=1}^nD_j^*D_j=L.
\]
Therefore
\begin{equation}\label{eq:first-form}
 \ip{Lv}{w}=\ip{Dv}{Dw}
     =\frac12\int\nabla v\cdot\nabla w\dga.
\end{equation}

To pass from this identity to the self-adjoint realization of $L$, let
$H=L^2(\ga)$ and define $\U=W^{1,2}(\ga)$ as the completion of
$C_c^\infty(\R^n)$ in the norm
\[
 \norm{v}_{\U}^2
   =\norm{v}_2^2+\frac12\norm{\nabla v}_2^2.
\]
Equivalently, $\U$ consists of the functions in $L^2(\ga)$ whose weak
Euclidean derivatives lie in $L^2(\ga)$.  For completeness, if $v$ and its weak derivatives are in $L^2(\ga)$,
choose $\chi_R\in C_c^\infty$ equal to one on $B_R$, supported in
$B_{2R}$, with $|\nabla\chi_R|\leq C/R$. Then
\[
 \norm{\nabla(\chi_Rv)-\nabla v}_2
 \leq\norm{(\chi_R-1)\nabla v}_2
       +\frac CR\norm{v\one_{B_{2R}\setminus B_R}}_2
 \longrightarrow0,
\]
and $\chi_Rv\to v$ in $L^2(\ga)$. On each fixed compact set the
Gaussian density is bounded above and below by positive constants, so
ordinary mollification finishes the approximation. Conversely,
convergence in this norm gives local unweighted $L^2$ convergence of
the functions and their derivatives, and hence the limiting weak
derivative identity. This also proves that
$D:\U\subset L^2(\ga)\to L^2(\ga;\R^n)$ is closed.

On $\U$ define the symmetric bilinear form
\[
 \mathfrak a(v,w)=\ip{Dv}{Dw}
   =\frac12\int\nabla v\cdot\nabla w\dga.
\]
It is densely defined because $\U$ is dense in $H$, nonnegative because
$\mathfrak a(v,v)=\norm{Dv}_2^2\geq0$, and closed because $\U$ is
complete in
$\bigl(\norm{v}_H^2+\mathfrak a(v,v)\bigr)^{1/2}$.  Together with
symmetry, these are exactly the hypotheses of the representation theorem
for closed forms; see \cite[Chapter~VI, \S2]{Kato95}. The theorem says that such an energy form determines
a unique nonnegative self-adjoint operator $A$ on $H$.  Concretely,
$v\in\U$ belongs to $\Dom(A)$ precisely when the functional
$w\mapsto\mathfrak a(v,w)$ can be represented by inner product with an
element $f\in H$; the operator then assigns $Av=f$.  Thus
\[
 \begin{split}
 \Dom(A)=\{v\in\U:\ &\text{there is an }f\in H\text{ such that}\\[-2pt]
                    &\mathfrak a(v,w)=\ip{f}{w}
                      \text{ for every }w\in\U\},
 \qquad Av=f.
 \end{split}
\]
Here nonnegative means $\ip{Av}{v}\geq0$, and self-adjoint means
$A=A^*$ with equal domains.  The representing function $f$ is unique
because $\U$ is dense in $H$. For a fixed $v\in C_c^\infty$,
identity \eqref{eq:first-form} extends from smooth test functions to
all $w\in\U$ by density. Consequently, every
$v\in C_c^\infty(\R^n)$ belongs to $\Dom(A)$ and
\[
 Av=-\frac12\Delta v+x\cdot\nabla v.
\]
Thus $A$ is the closed self-adjoint realization denoted by $L$.

The square-root part of the representation theorem further gives
\[
 \Dom(L^{1/2})=\U,
 \qquad
 \mathfrak a(v,w)=\ip{L^{1/2}v}{L^{1/2}w}.
\]
Taking $w=v$ yields the fundamental identity
\begin{equation}\label{eq:sqrt-domain}
 \begin{split}
 \U&=\Dom(L^{1/2})=W^{1,2}(\ga),\\
 \norm{L^{1/2}v}_2^2&=\norm{Dv}_2^2
   =\frac12\int|\nabla v|^2\dga.
 \end{split}
\end{equation}
Thus $\Dom(L^{1/2})$ requires one weak derivative in $L^2(\ga)$,
whereas $\Dom(L)$ additionally requires
$-\frac12\Delta v+x\cdot\nabla v\in L^2(\ga)$.  Constants belong to the
square-root domain and are annihilated by both $L^{1/2}$ and $D$.
The Hermite expansion below provides the equivalent spectral description
\[
 v\in\Dom(L^{1/2})
 \quad\Longleftrightarrow\quad
 \sum_{\alpha\ne0}|\alpha||v_\alpha|^2<\infty.
\]
We reserve $\U$ for this first-order space; $\V$ will denote the distinct
fractional form domain introduced in \cref{subsec:fractional-domain}.
Standard Gaussian operator conventions may also be found in
\cite{FGS94,Bruno19}.

\subsection{Hermite expansion and constants}
\label{subsec:hermite}

For a multi-index $\alpha\in\N_0^n$, where
$\N_0=\{0,1,2,\ldots\}$, write
$|\alpha|=\alpha_1+\cdots+\alpha_n$.
Let $H_k(t)=(-1)^k e^{t^2}\frac{d^k}{dt^k}e^{-t^2}$ denote the
one-dimensional Hermite polynomials for the weight $e^{-t^2}$.
Throughout the paper, the \emph{Hermite basis} means the normalized
multidimensional Hermite polynomials, defined as the tensor products
\[
 h_\alpha(x)=\prod_{j=1}^n
     \frac{H_{\alpha_j}(x_j)}{\sqrt{2^{\alpha_j}\alpha_j!}},
 \qquad x=(x_1,\ldots,x_n)\in\R^n.
\]
These polynomials form an orthonormal basis of $L^2(\ga)$; see
\cite{Urbina19}. Write also $H_\alpha=\prod_{j=1}^nH_{\alpha_j}$
for the unnormalized tensor polynomial. We have $h_0=1$ and
\begin{equation}\label{eq:Hermite}
 Lh_\alpha=|\alpha|h_\alpha.
\end{equation}
For this normalization,
$D_jh_\alpha=\sqrt{\alpha_j}\,h_{\alpha-e_j}$, where $e_j$ is the
$j$th coordinate vector and the expression is zero if $\alpha_j=0$.

We verify that these are eigenfunctions of the closed-form realization,
not only of the formal differential expression. Every polynomial and
its first derivatives belong to $L^2(\ga)$, hence $h_\alpha\in\U$.
Integration by parts first against $w\in C_c^\infty$, and then density
in $\U$, give
\[
 \mathfrak a(h_\alpha,w)=|\alpha|\ip{h_\alpha}{w},
 \qquad w\in\U.
\]
The representation of the form therefore implies
$h_\alpha\in\Dom(L)$ and \eqref{eq:Hermite}. The diagonal operator
with eigenvalues $|\alpha|$ in this complete orthonormal basis is
self-adjoint. If $v\in\Dom(L)$, self-adjointness gives
$\ip{Lv}{h_\alpha}=|\alpha|\ip{v}{h_\alpha}$, so $L$ is contained
in that diagonal operator. Two self-adjoint operators related by
extension coincide. This identifies the closed-form realization with
the Hermite spectral realization.

Writing $v_\alpha=\ip{v}{h_\alpha}$, spectral calculus gives, for $s>0$,
\begin{equation}\label{eq:spectral-domains}
 \begin{split}
 \Dom(L^s)
   &=\left\{v\in L^2(\ga):
       \sum_{\alpha\ne0}|\alpha|^{2s}|v_\alpha|^2<\infty\right\},\\
 L^sv&=\sum_{\alpha\ne0}|\alpha|^sv_\alpha h_\alpha.
 \end{split}
\end{equation}
Thus constants belong to every positive-power domain and are annihilated
by $L^s$. To separate the mean-zero part from constants, define the mean
and the mean-zero projection by
\begin{equation}\label{eq:projection}
 \me(v)=\int v\dga,
 \qquad \Pi_0v=v-\me(v)1.
\end{equation}
Since the constant eigenvalue is zero, we define the inverse square root
only on mean-zero functions and set it to zero on constants:
\begin{equation}\label{eq:pseudoinverse}
 L^{-1/2}\Pi_0v
   =\sum_{\alpha\ne0}|\alpha|^{-1/2}v_\alpha h_\alpha.
\end{equation}
For every $v\in L^2(\ga)$, this series converges in $L^2(\ga)$ because
$|\alpha|\geq1$ for $\alpha\ne0$; its norm is at most
$\norm{\Pi_0v}_2$. Moreover, the potential $w=L^{-1/2}\Pi_0v$ belongs
to $\Dom(L^{1/2})$, since
\[
 \sum_{\alpha\ne0}|\alpha||w_\alpha|^2
 =\sum_{\alpha\ne0}|v_\alpha|^2
 =\norm{\Pi_0v}_2^2<\infty.
\]
Applying \eqref{eq:spectral-domains} now gives
\[
 L^{1/2}w=\Pi_0v,
 \qquad \norm{L^{1/2}w}_2=\norm{\Pi_0v}_2.
\]
In other words, $w$ is the unique mean-zero solution of
$L^{1/2}w=v-\me(v)1$.

\begin{definition}[Energy-normalized vector Riesz transform]
\label{def:riesz}
For $f\in L^2(\ga)$, set
\begin{equation}\label{eq:riesz-def}
 \Rg f=D L^{-1/2}\Pi_0f.
\end{equation}
Its $j$th coordinate is $R_{\gamma,j}f=D_jL^{-1/2}\Pi_0f$.
The preceding calculation and \eqref{eq:sqrt-domain} show that the
potential $L^{-1/2}\Pi_0f$ lies in $W^{1,2}(\ga)$, so its weak gradient
belongs to $L^2(\ga;\R^n)$. Thus we differentiate the potential, not
$f$ itself; no differentiability assumption on the input is needed.
\end{definition}

\begin{proposition}[The two operator identities]
\label{prop:identities}
For every $f\in L^2(\ga)$,
\begin{equation}\label{eq:R-L2}
 \norm{\Rg f}_{L^2(\ga;\R^n)}^2
    =\norm{\Pi_0f}_2^2\leq\norm{f}_2^2.
\end{equation}
For every $u\in\Dom(L^{1/2})$,
\begin{equation}\label{eq:R-B}
 \Rg L^{1/2}u=Du
 \quad\text{in }L^2(\ga;\R^n).
\end{equation}
\end{proposition}

\begin{proof}
For $v=L^{-1/2}\Pi_0f$, we have already shown that
$v\in\Dom(L^{1/2})$ and $L^{1/2}v=\Pi_0f$. The energy identity
\eqref{eq:sqrt-domain} therefore gives
\[
 \norm{\Rg f}_2^2=\norm{Dv}_2^2
       =\norm{L^{1/2}v}_2^2=\norm{\Pi_0f}_2^2.
\]
For $u\in\Dom(L^{1/2})$, the Hermite expansion gives
\[
 L^{-1/2}\Pi_0L^{1/2}u=\Pi_0u.
\]
Applying $D$ gives $D\Pi_0u=Du$, since $D1=0$.
\end{proof}

\begin{remark}[The projection and the input norm]
The projection in \eqref{eq:riesz-def} is essential, but we will decompose
$f$ itself, using its nonnegative parts $f^\pm$, not center it first.
This avoids the unnecessary estimate $\norm{\Pi_0f}_1\leq2\norm{f}_1$;
the identities above already handle constants.
\end{remark}

\subsection{Relation with the normalization in the 1994 paper}
\label{subsec:FGS-normalization}

On page 232 of \cite{FGS94}, the spectral convention, for
$|\alpha|>0$, is
\[
 \mathcal R_jH_\alpha
    =-\frac1{\sqrt{|\alpha|}}\partial_jH_\alpha.
\]
Thus, when both operators are set to zero on constants,
\begin{equation}\label{eq:FGS-spectral}
 \mathcal R_j=-\partial_jL^{-1/2}\Pi_0=-\sqrt2\,R_{\gamma,j}.
\end{equation}
The sign of the Ornstein--Uhlenbeck generator in that paper is the
opposite of our nonnegative $L$. The factor $\sqrt2$ in
\eqref{eq:FGS-spectral} comes from $\partial_j=\sqrt2 D_j$, not
from this sign convention. Consequently, \eqref{eq:target} gives
$2\sqrt2$ for its vector spectral normalization.
Multiplying the underlying measure by $\pi^{n/2}$ to use
$e^{-|x|^2}\,dx$ does not change a weak-type constant for a fixed operator:
both the distribution measure and the input $L^1$ norm acquire the same
factor.

\paragraph{The normalized kernel.}
The \emph{Mehler semigroup}, also called the Ornstein--Uhlenbeck
semigroup, is the family $P_s=e^{-sL}$, $s\geq0$. For
$f\in L^2(\ga)$ it is defined by the Hermite expansion
\[
 P_sf=\sum_{\alpha\in\N_0^n}e^{-s|\alpha|}f_\alpha h_\alpha,
 \qquad f_\alpha=\ip{f}{h_\alpha},
\]
with convergence in $L^2(\ga)$. Thus it preserves the constant Hermite
component and damps the nonconstant Hermite components. Moreover,
$P_0=I$ and $P_sP_t=P_{s+t}$ for $s,t\geq0$. Strong continuity at
zero follows from Parseval and dominated convergence applied to
$|e^{-s|\alpha|}-1|^2|f_\alpha|^2$.
Its explicit averaging formula, also used in \cref{subsec:Markov} below, is
\[
 P_sf(x)=\int_{\R^n}
 f\bigl(e^{-s}x+\sqrt{1-e^{-2s}}\,z\bigr)\dga(z),
 \qquad s>0.
\]
This is Mehler's formula; see \cite{Urbina19}. Its identification with
the spectral semigroup can also be checked from the generating function
$\sum_{k\geq0}H_k(t)z^k/k!=e^{2tz-z^2}$. Gaussian integration
shows that averaging $H_k(rx+\sqrt{1-r^2}\,y)$ gives $r^kH_k(x)$.
Tensorization gives the same action on every $h_\alpha$. The averaging
operator is an $L^2$ contraction by Jensen's inequality and Gaussian
invariance, so equality on the Hermite span proves equality on $L^2$.
In words, it averages $f$ after contracting $x$ toward the origin
and adding Gaussian noise. By the \emph{Lebesgue kernel} we mean
the function $p_s$ such that
\[
 P_sf(x)=\int_{\R^n}p_s(x,y)f(y)\,dy.
\]
Here the integration is with respect to Lebesgue measure $dy$, not
$d\ga(y)$. Changing variables in the averaging formula gives
\begin{equation}\label{eq:Mehler-Lebesgue}
 p_s(x,y)=\frac{1}{\pi^{n/2}(1-e^{-2s})^{n/2}}
   \exp\left(-\frac{|e^{-s}x-y|^2}{1-e^{-2s}}\right).
\end{equation}
For every $f\in L^2(\ga)$, spectral calculus gives the convergent
Bochner integral
\begin{equation}\label{eq:subordination-centered}
 L^{-1/2}\Pi_0f
 =\frac1{\sqrt\pi}\int_0^\infty
       (P_sf-\me(f)1)\,\frac{ds}{\sqrt s}.
\end{equation}
Indeed, $P_sf-\me(f)1=P_s\Pi_0f$ has $L^2$ norm at most
$\norm{f}_2$ for $0<s\leq1$ and at most
$e^{-s}\norm{\Pi_0f}_2$ for $s\geq1$. Without this subtraction,
the integral diverges at infinity when $\me(f)\ne0$.
For $f\in C_c^\infty$, differentiation of Mehler's formula gives
$DP_sf=e^{-s}P_sDf$. Truncating the integral in
\eqref{eq:subordination-centered}, and then using closedness of $D$
and $L^2$ convergence of $\int e^{-s}P_sDf\,ds/\sqrt{s}$, yields
\begin{equation}\label{eq:R-gradient-semigroup}
 \Rg f=\frac1{\sqrt\pi}\int_0^\infty
      DP_sf\,\frac{ds}{\sqrt s}
 =\frac1{\sqrt\pi}\int_0^\infty
      e^{-s}P_sDf\,\frac{ds}{\sqrt s}.
\end{equation}
The second integral converges absolutely at every $x$, since
$|P_sDf(x)|\leq\norm{Df}_\infty$, and provides a continuous
representative of $\Rg f$.

Differentiating the semigroup, or first truncating the time integral,
and setting $r=e^{-s}$ shows that the off-diagonal Lebesgue kernel of
$R_{\gamma,j}$ is
\begin{equation}\label{eq:normalized-kernel}
 \begin{split}
 K_{\gamma,j}(x,y)
 &=\frac{\sqrt2}{\pi^{(n+1)/2}}\,k_j^{\rm raw}(x,y),\\
 k_j^{\rm raw}(x,y)
 &=\int_0^1
   \frac{y_j-rx_j}{(1-r^2)^{n/2+1}\sqrt{-\log r}}
   \exp\left(-\frac{|rx-y|^2}{1-r^2}\right)\,dr,
   \qquad x\ne y.
 \end{split}
\end{equation}
Indeed,
$\partial_{x_j}p_s=2r(y_j-rx_j)(1-r^2)^{-1}p_s$;
the factor $r$ cancels the $1/r$ in $ds=-dr/r$.
For compactly supported smooth inputs and $x$ off their support, the
integrals can be interchanged directly. The formula is an off-diagonal
kernel identification, not yet a principal-value identification. The latter is proved in
\cref{prop:pv-identification,lem:pv-core}.

The expression $k_j^{\rm raw}$ is the kernel displayed at the beginning
of \cite{FGS94}, after combining its two powers of $1-r^2$.
That paper explicitly distinguishes its kernel and spectral definitions
up to a multiplicative constant. Formula \eqref{eq:normalized-kernel}
restores the prefactor needed here. In particular, the kernel of
$\mathcal R_j$ in \eqref{eq:FGS-spectral} is
$-2\pi^{-(n+1)/2}k_j^{\rm raw}$. A statement about dimension-free constants
must specify these factors; an unnormalized kernel is a different family
of operators.

\subsection{The standard-variance Gaussian convention}

Under the change of variables $z=\sqrt2\,x$, define
$Uf(z)=f(z/\sqrt2)$. Then $U$ preserves all $L^p$ norms and distribution
functions when $\ga$ is replaced by
$(2\pi)^{-n/2}e^{-|z|^2/2}\,dz$. Moreover,
\[
 UL=(-\Delta_z+z\cdot\nabla_z)U,
 \qquad UD=\nabla_zU.
\]
Thus the constant $2$ corresponds to
$\nabla_z(-\Delta_z+z\cdot\nabla_z)^{-1/2}\Pi_0$ in standard-variance
coordinates. The normalization is not an additional dimensional loss.

\section{The Gaussian obstacle decomposition and the height restriction}
\label{sec:decomposition}

Write
\begin{equation}\label{eq:B-def}
 B=L^{1/2}.
\end{equation}
The following is the exact, unshifted Gaussian obstacle decomposition. Its proof occupies \cref{sec:proof}.

\begin{theorem}[Gaussian obstacle decomposition]
\label{thm:decomposition}
Let $F\in L^2(\ga)$ be nonnegative, put $m=\int F\dga$, and let $a>m$.
There exist $u\in\Dom(B)=W^{1,2}(\ga)$ and
$g\in L^2(\ga)\cap L^\infty(\ga)$ such that
\begin{align}
 &u\geq0,\qquad F=g+Bu \quad\text{in }L^2(\ga),
       \label{eq:decomp-equality}\\
 &\min\{F,a\}\leq g\leq a \quad\text{a.e.},
       \label{eq:g-bounds}\\
 &\int g\dga=m.
       \label{eq:g-mass}
\end{align}
For a measurable representative of $u$, set $\Omega=\{u>0\}$. Then
\begin{align}
 &g=a\quad\text{a.e. on }\Omega,
       \label{eq:g-contact}\\
 &a\,\ga(\Omega)\leq m,
       \label{eq:Omega-bound}\\
 &Du=0\quad\text{a.e. on }\Omega^c.
       \label{eq:gradient-support}
\end{align}
The pair $(u,g)$ satisfying these assertions is unique. The potential
$u$ is the unique minimizer of the functional in \eqref{eq:functional}
over the cone in \eqref{eq:cone}.
\end{theorem}

The bad function $b=F-g=Bu$ has zero Gaussian mean. It is \emph{not}
claimed to be supported on $\Omega$. Instead, the property relevant to
the endpoint estimate is
\begin{equation}\label{eq:bad-transform}
 \Rg b=Du=0\quad\text{a.e. on }\Omega^c.
\end{equation}
This distinction is the reason the decomposition is adapted to the Riesz
transform.

\subsection{Why the upper bound $a$ on $g$ cannot be arbitrary}

Here $a$ is the upper bound prescribed for $g$, and $m=\int F\dga$.
The assumption $a>m$ in \cref{thm:decomposition} is essential: if
$a<m$, the bound $g\leq a$ is incompatible with mass preservation,
while at $a=m$ any minimizer,
if it exists, is not unique. These obstructions come from the finite
total mass of $\ga$ and the fact that $B=L^{1/2}$ annihilates constants.
Indeed, if
$F=g+Bu$ with $u\in\Dom(B)$, then self-adjointness of $B$ in
$L^2(\ga)$ and $1\in\Dom(B)$ give
\[
 \int Bu\dga=\ip{Bu}{1}=\ip{u}{B1}=0.
\]
So any such decomposition automatically preserves mass: $\int g\dga=\int F\dga$.

If in addition $0\leq g\leq a$, then, because $\ga$ has mass one,
\[
 m=\int g\dga\leq\int a\dga=a.
\]
Thus no mass-preserving decomposition with this upper bound can exist
when $a<m$. For example, if $F\equiv m>a$, no function bounded above
by $a$ can have the same Gaussian integral as $F$. This restriction
uses $\ga(\R^n)=1$: on an infinite-measure space the bound $g\leq a$
alone would not give a finite upper bound on $\int g$.

The strict inequality $a>m$, rather than merely $a\geq m$, also
controls the constant part of the obstacle potential. To see this,
consider the functional $J$ in \eqref{eq:functional}, defined on
nonnegative elements of $\V=\Dom(L^{1/4})$. Its quadratic energy
$\E(v)=\norm{L^{1/4}v}_2^2$ is unchanged by adding a constant,
because $L^{1/4}1=0$. Hence, for every admissible $v$ and constant
$c\geq0$, the function $v+c$ is still admissible and
\[
 J(v+c)=J(v)+(a-m)c.
\]
In particular, the constant competitors satisfy $J(c)=(a-m)c$.
If $a<m$, these values tend to $-\infty$ as $c\to\infty$, so there
can be no minimizer. If $a>m$, adding a large positive constant has
a positive cost. Together with the spectral-gap control of the
mean-zero part, this yields the coercivity proved in
\cref{subsec:existence}.

At the borderline height $a=m$, the mass constraint and $g\leq a$
would force $g=a$ almost everywhere: the nonnegative function $a-g$
has integral $a-m=0$. The decomposition would therefore require
$Bu=F-m$. Moreover, $J(v+c)=J(v)$, so the functional does not
control constants; if a minimizer exists, adding any positive constant
produces another one. This explains why the existence and uniqueness
argument below assumes $a>m$. We make no general assertion of existence
at $a=m$.

Nor should one impose zero mean on the nonnegative obstacle potential:
$u\geq0$ and $\int u\dga=0$ would force $u=0$ almost everywhere.
Indeed, for every $\varepsilon>0$,
\[
 0=\int u\dga\geq\varepsilon\,\ga(\{u>\varepsilon\}),
\]
and $\{u>0\}=\bigcup_{k\geq1}\{u>1/k\}$ is therefore null.
The decomposition would then reduce to $g=F$, which cannot satisfy
$g\leq a$ when $F>a$ on a set of positive measure. Subtracting the
mean from a nonconstant nonnegative potential would therefore violate
the nonnegativity constraint. The minimization instead keeps that
constant part of $u$ and controls it by the positive quantity $a-m$.

\subsection{Why the height restriction does not obstruct weak type}

For a real $f\in L^2(\ga)$, put $M=\norm{f}_1$. If $M=0$, then
$f=0$ almost everywhere. For $0<t\leq M$, the estimate is immediate
because $\ga$ is a probability measure:
\begin{equation}\label{eq:small-level}
 t\,\ga(\{|\Rg f|>t\})\leq t\leq M.
\end{equation}
In particular, the borderline level $t=M$ needs no decomposition.
For $t>M$, both $f^+$ and $f^-$ have mean at most $M<t$, so
\cref{thm:decomposition} applies to both at the common height $a=t$.
This is a condition on their integrals, not on their pointwise values:
the inputs can be unbounded. The two decompositions are combined
explicitly in \cref{subsec:signed}. Once the estimate has been proved
on $L^2$, \cref{subsec:L1extension} constructs the $L^1$ extension.

\section{Detailed proof of the Gaussian obstacle decomposition}
\label{sec:proof}

\paragraph{Plan of the proof.}
For $F\geq0$ in $L^2(\ga)$ and $a>m=\int F\dga$, we construct
$u$ and $g$ and verify all assertions of \cref{thm:decomposition}
in six steps.
\begin{enumerate}[label=\textup{(\arabic*)},leftmargin=*,itemsep=2pt]
 \item Establish the fractional form domain $\V=\Dom(L^{1/4})$,
 its spectral gap, and its Markov and lattice properties.
 \item Obtain the unique minimizer $u\in\K\subset\V$ of $J$.
 The spectral gap controls the mean-zero part, and $a>m$ controls
 the constant part.
 \item Derive the variational inequality, introduce the multiplier
 $\eta=J'(u)\in\V^*$, and prove the Lewy--Stampacchia bound
 $0\leq\eta\leq(a-F)^+$ in the dual order
 \eqref{eq:eta-order-V}.
 \item Represent $\eta$ in $L^2$, define $g=a-\eta$, and use the
 form equation to upgrade $u$ to $\U=\Dom(B)=W^{1,2}(\ga)$.
 This gives $F=g+Bu$ and $\min\{F,a\}\leq g\leq a$ a.e.
 \item Test against constants to prove $\int g\dga=m$, and obtain
 complementarity $\eta u=0$ a.e. For $\Omega=\{u>0\}$, this yields
 $g=a$ a.e. on $\Omega$ and $a\ga(\Omega)\leq m$.
 \item Use the Sobolev chain rule and closedness of $D$ to prove
 $Du=0$ a.e. on $\Omega^c$. Together with \eqref{eq:R-B}, this gives
 \[
  \Rg(F-g)=Du=\frac1{\sqrt2}\nabla u=0
  \quad\text{a.e. on }\Omega^c,
 \]
 the localization needed for the endpoint estimate in
 \cref{sec:endpoint}.
\end{enumerate}
The domain upgrade in \cref{subsec:upgrade} must precede the use of $Bu$, $Du$,
and the Riesz-transform identity: minimization initially gives only
$u\in\V$.

\subsection{The correct fractional form domain}
\label{subsec:fractional-domain}

The operator $B=L^{1/2}$ has operator domain $\Dom(B)=\U$: these are
the functions $v\in L^2(\ga)$ for which $Bv\in L^2(\ga)$.
Its \emph{form domain}, by contrast, consists of the functions
$v\in L^2(\ga)$ with finite associated energy
$\norm{B^{1/2}v}_2^2$. Since $B^{1/2}=L^{1/4}$, this larger space
and the corresponding bilinear form are
\begin{equation}\label{eq:V-form}
 \V=\Dom(L^{1/4}),\qquad
 \E(v,w)=\ip{L^{1/4}v}{L^{1/4}w},\qquad
 \E(v)=\E(v,v).
\end{equation}
Here ``form'' refers to the scalar-valued bilinear expression $\E(v,w)$,
and $\E(v)$ is its energy. Having finite energy does not by itself
guarantee that $Bv$ is an $L^2$ function.

Equip $\V$ with the Hilbert norm
$\norm{v}_{\V}^2=\norm{v}_2^2+\E(v)$. In Hermite coordinates,
\begin{equation}\label{eq:form-series}
 \E(v)=\sum_{\alpha\ne0}|\alpha|^{1/2}|v_\alpha|^2.
\end{equation}
Here $v_\alpha=\ip{v}{h_\alpha}=\int v h_\alpha\dga$ is the
coefficient of $v$ in the normalized multidimensional Hermite basis of
$L^2(\ga)$ introduced in \cref{subsec:hermite}.
The index $\alpha=0$ corresponds to $h_0=1$ and
$v_0=\me(v)$; it contributes no energy and is omitted from the sum.
The two spaces are related by
\[
 \U\subsetneq\V\subset L^2(\ga),
 \qquad \norm{v}_{\V}\leq\norm{v}_{\U}\quad(v\in\U).
\]
The inclusion and norm bound follow from
$|\alpha|^{1/2}\leq|\alpha|$ for $\alpha\ne0$; strictness is verified in
\cref{rem:quarter-power-domain}.  Thus $\U$ is the form domain of $L$,
whereas $\V$ is the form domain of $B=L^{1/2}$, not another name for $\U$.

\begin{remark}[Meaning of the quarter-power domain]
\label{rem:quarter-power-domain}
The meaning of $\Dom(L^{1/4})$ is given by the spectral-calculus
criterion \eqref{eq:spectral-domains} with $s=1/4$; in particular,
\eqref{eq:form-series} says exactly that
$\E(v)=\norm{L^{1/4}v}_2^2$. Thus the $\V$-norm introduced above has
the Hermite-coordinate expression
\[
 \norm{v}_{\V}^2=|v_0|^2+\sum_{\alpha\ne0}
       \bigl(1+|\alpha|^{1/2}\bigr)|v_\alpha|^2,
\]
and its inner product is $\langle v,w\rangle_{\V}
=\ip{v}{w}+\E(v,w)$. The weighted $\ell^2$ description shows that this
norm makes $\V$ a Hilbert space. Its energy alone is only a seminorm:
the constant Hermite component has zero energy, but is controlled by the
$L^2$ term in the $\V$-norm.

This quarter-power domain is the form domain of $B=L^{1/2}$, since the
form domain of a nonnegative self-adjoint operator $B$ is
$\Dom(B^{1/2})$. Although $\E(v,w)=\ip{Bv}{w}$ when
$v\in\Dom(B)=\U$ and $w\in\V$, that $L^2$-pairing is not available for
arbitrary $v\in\V$.

The inclusion $\U\subsetneq\V$ is strict. Indeed, for any dimension $n$
let $v=\sum_{k\ge1}k^{-7/8}h_{(k,0,\ldots,0)}$. Its coefficient
squares are summable, and the domain sums in \eqref{eq:spectral-domains}
for $s=1/4$ and $s=1/2$ are respectively
$\sum_{k\ge1}k^{-5/4}<\infty$ and
$\sum_{k\ge1}k^{-3/4}=\infty$. Hence $v\in\V\setminus\U$; by
\eqref{eq:sqrt-domain}, it has no $L^2$ weak gradient. This distinction
matters in the obstacle proof: minimization initially gives only
$u\in\V$, and the upgrade to $u\in\U$ in \cref{subsec:upgrade} is
needed before using $Bu$ or the gradient identity as $L^2$ quantities.
\end{remark}

A direct consequence of \eqref{eq:form-series} is the fractional
spectral-gap estimate
\begin{equation}\label{eq:fractional-gap}
 \norm{v-\me(v)}_2^2\leq\E(v),\qquad v\in\V.
\end{equation}
Indeed, every nonconstant Hermite basis function $h_\alpha$ has
eigenvalue $|\alpha|\geq1$, so $|\alpha|^{1/2}\geq1$.
Also, $1\in\V$ and $\E(v,1)=0$ for every $v\in\V$.

\subsection{Mehler representation, positivity, and Gaussian invariance}
\label{subsec:Markov}

For $s>0$, the Ornstein--Uhlenbeck semigroup has the probability formula
\begin{equation}\label{eq:Mehler-prob}
 P_sv(x)=\int v\bigl(e^{-s}x+\sqrt{1-e^{-2s}}\,y\bigr)\dga(y).
\end{equation}
In particular, it is positivity preserving, preserves constants, and
preserves Gaussian integrals. The last assertion means that, for every
$v\in L^2(\ga)$ and $s>0$,
\[
 \int_{\R^n}P_sv(x)\dga(x)
 =\int_{\R^n}v(x)\dga(x),
 \qquad\text{equivalently, }\me(P_sv)=\me(v).
\]
Thus applying $P_s$ leaves the Gaussian mean of $v$ unchanged, although
it generally changes $v$ itself. To verify this, the Hermite expansion
shows that $P_s$ is self-adjoint on $L^2(\ga)$ and $P_s1=1$, so
\[
 \int P_sv\dga=\ip{P_sv}{1}=\ip{v}{P_s1}
 =\ip{v}{1}=\int v\dga.
\]
This property is also called invariance of the measure $\ga$ under
$P_s$.
Relative to $\ga$, the semigroup is represented explicitly by
\[
 P_sv(x)=\int_{\R^n}M_s(x,y)v(y)\dga(y),
 \qquad s>0,\quad v\in L^2(\ga),
\]
where its Gaussian-measure kernel is
\begin{equation}\label{eq:Mehler-Gaussian}
 M_s(x,y)=\frac{1}{(1-r^2)^{n/2}}
   \exp\left(\frac{2r\,x\cdot y-r^2(|x|^2+|y|^2)}{1-r^2}\right),
 \qquad r=e^{-s}.
\end{equation}
It is positive and symmetric, and its integral in either variable is
one. Jensen's inequality and invariance of $\ga$ give $L^2$ contraction.
The spectral action is $P_sh_\alpha=e^{-s|\alpha|}h_\alpha$.

For $\lambda\geq0$, integration by parts gives the scalar identity
\begin{equation}\label{eq:scalar-fractional}
 \sqrt\lambda=\frac1{2\sqrt\pi}
    \int_0^\infty(1-e^{-s\lambda})\,\frac{ds}{s^{3/2}}.
\end{equation}
For $\lambda>0$, the integral equals
$2\lambda\int_0^\infty e^{-s\lambda}s^{-1/2}\,ds
=2\sqrt\pi\sqrt\lambda$; the case $\lambda=0$ is obvious.

We now apply formula \eqref{eq:scalar-fractional} to each eigenvalue
in the spectral expansion. For $v\in L^2(\ga)$, write
$v=\sum_\alpha v_\alpha h_\alpha$, where
$v_\alpha=\ip{v}{h_\alpha}$ and $\sum_\alpha|v_\alpha|^2=\norm{v}_2^2$.
Since $P_sh_\alpha=e^{-s|\alpha|}h_\alpha$, orthonormality gives,
for each $s>0$,
\[
 \ip{(I-P_s)v}{v}
 =\sum_{\alpha\ne0}(1-e^{-s|\alpha|})|v_\alpha|^2.
\]
The constant component contributes zero because $P_s1=1$, so
$(I-P_s)1=0$. This series converges, since
$0\leq1-e^{-s|\alpha|}\leq1$, and its sum lies between $0$ and
$\norm{v}_2^2$. Each summand is a nonnegative measurable function of
$s$. Tonelli's theorem, applied to counting measure on the Hermite
indices and the measure $s^{-3/2}\,ds$ on $(0,\infty)$, therefore
allows us to interchange the sum and the integral without first
assuming that the integral is finite:
\[
 \begin{aligned}
 \frac1{2\sqrt\pi}\int_0^\infty
       \ip{(I-P_s)v}{v}\,\frac{ds}{s^{3/2}}
 &=\sum_{\alpha\ne0}|v_\alpha|^2
       \frac1{2\sqrt\pi}\int_0^\infty
       (1-e^{-s|\alpha|})\,\frac{ds}{s^{3/2}}\\
 &=\sum_{\alpha\ne0}|\alpha|^{1/2}|v_\alpha|^2.
 \end{aligned}
\]
In the last equality we used \eqref{eq:scalar-fractional} with
$\lambda=|\alpha|$. By \eqref{eq:form-series}, the resulting sum is
$\E(v)$ when $v\in\V=\Dom(L^{1/4})$; by the spectral domain
criterion, it is $+\infty$ when $v\in L^2(\ga)\setminus\V$.
With the convention $\E(v)=+\infty$ in the latter case, we obtain
\begin{equation}\label{eq:semigroup-energy}
 \E(v)=\frac1{2\sqrt\pi}\int_0^\infty
          \ip{(I-P_s)v}{v}\,\frac{ds}{s^{3/2}},
 \qquad v\in L^2(\ga).
\end{equation}
For clarity, symmetry means $M_s(x,y)=M_s(y,x)$, and the two
unit-mass identities are
\[
 \int_{\R^n}M_s(x,y)\dga(y)=1\quad\text{for each }x,
 \qquad
 \int_{\R^n}M_s(x,y)\dga(x)=1\quad\text{for each }y.
\]
Both identities follow from $P_s1=1$ and the symmetry of $M_s$.
For real $v\in L^2(\ga)$, Tonelli's theorem gives
\[
 \begin{aligned}
 \iint v(x)^2M_s(x,y)\dga(x)\dga(y)
 &=\int v(x)^2\left(\int M_s(x,y)\dga(y)\right)\dga(x)
 =\norm{v}_2^2,\\
 \iint v(y)^2M_s(x,y)\dga(x)\dga(y)
 &=\int v(y)^2\left(\int M_s(x,y)\dga(x)\right)\dga(y)
 =\norm{v}_2^2.
 \end{aligned}
\]
The mixed term is absolutely integrable, since
$2|v(x)v(y)|\leq v(x)^2+v(y)^2$ and the preceding two integrals
are finite. Fubini's theorem and the kernel representation of $P_s$
therefore yield
\[
 \iint v(x)v(y)M_s(x,y)\dga(x)\dga(y)
 =\int v(x)P_sv(x)\dga(x)=\ip{P_sv}{v}.
\]
Expanding $|v(x)-v(y)|^2=v(x)^2+v(y)^2-2v(x)v(y)$ now gives
\[
 \begin{aligned}
 \frac12\iint |v(x)-v(y)|^2M_s(x,y)\dga(x)\dga(y)
 &=\frac12\bigl(2\norm{v}_2^2-2\ip{P_sv}{v}\bigr)\\
 &=\ip{(I-P_s)v}{v}.
 \end{aligned}
\]
Substituting this identity into \eqref{eq:semigroup-energy} gives
\begin{equation}\label{eq:difference-energy}
 \E(v)=\frac1{4\sqrt\pi}
 \int_0^\infty\iint |v(x)-v(y)|^2M_s(x,y)
             \dga(x)\dga(y)\,\frac{ds}{s^{3/2}}.
\end{equation}
By polarization and Cauchy--Schwarz for the measure appearing here,
for $v,w\in\V$ one also has
\begin{equation}\label{eq:bilinear-energy}
 \E(v,w)=\frac1{2\sqrt\pi}\int_0^\infty
                 \ip{(I-P_s)v}{w}\,\frac{ds}{s^{3/2}}.
\end{equation}
The bilinear integral is absolutely convergent. The role of the kernel
here differs from the local--global Riesz-kernel analysis in
\cite[pp.~230--233]{FGS94}: positivity, symmetry, and unit mass give
\eqref{eq:difference-energy}, and hence the lattice properties below,
without a pointwise size estimate for $M_s$.

\begin{lemma}[Lattice properties and the cross-term sign]
\label{lem:lattice}
If $\Phi:\R\to\R$ is $1$-Lipschitz with $\Phi(0)=0$ and $v\in\V$,
then $\Phi(v)\in\V$ and $\E(\Phi(v))\leq\E(v)$.
In particular, $v^+$, $v^-$, and $|v|$ belong to $\V$.
If $v,w\in\V$ are nonnegative and $vw=0$ a.e., then
$\E(v,w)\leq0$. Consequently,
\begin{equation}\label{eq:negative-part-sign}
 \E(z,z^-)\leq-\E(z^-),\qquad z\in\V,
\end{equation}
where $z^-=\max\{-z,0\}$.
\end{lemma}

\begin{proof}
Since $|\Phi(v)|\leq|v|$, the composition belongs to $L^2$.
The pointwise inequality
$|\Phi(v(x))-\Phi(v(y))|\leq|v(x)-v(y)|$ holds for almost every
$(x,y)$. Because $M_s(x,y)\geq0$, squaring this inequality,
multiplying by $M_s(x,y)$, and integrating first in $x,y$ and then in
$s$ gives
\[
 \begin{aligned}
 &\frac1{4\sqrt\pi}\int_0^\infty\iint
   |\Phi(v(x))-\Phi(v(y))|^2M_s(x,y)
   \dga(x)\dga(y)\,\frac{ds}{s^{3/2}}\\
 &\qquad\leq
 \frac1{4\sqrt\pi}\int_0^\infty\iint
   |v(x)-v(y)|^2M_s(x,y)
   \dga(x)\dga(y)\,\frac{ds}{s^{3/2}}
 =\E(v)<\infty.
 \end{aligned}
\]
Formula \eqref{eq:difference-energy} applies to every $L^2(\ga)$
function with the convention that its right-hand side is $+\infty$
outside $\V$. Since the displayed integral for $\Phi(v)$ is finite,
it follows that $\Phi(v)\in\V$; applying the same formula once more
identifies that integral with $\E(\Phi(v))$ and yields
$\E(\Phi(v))\leq\E(v)$. Taking respectively
$\Phi(t)=t^+$, $\Phi(t)=t^-$, and $\Phi(t)=|t|$ gives the stated
lattice properties, since all three maps are $1$-Lipschitz and vanish
at zero.
For nonnegative disjoint $v,w$, we have
\[
 \ip{(I-P_s)v}{w}=-\ip{P_sv}{w}\leq0.
\]
Integrating in \eqref{eq:bilinear-energy} gives $\E(v,w)\leq0$.
Finally, write $z=z^+-z^-$ and use
$\E(z^+,z^-)\leq0$.
\end{proof}

\subsection{Coercivity, existence, and uniqueness of a minimizer}
\label{subsec:existence}

\begin{lemma}[Coercivity and the unique obstacle minimizer]
\label{lem:obstacle-minimizer}
Let $F\in L^2(\ga)$ be nonnegative, let $m=\int F\dga$, and assume
$a>m$. Recall that $\V=\Dom(L^{1/4})$ with
$\norm{v}_{\V}^2=\norm{v}_2^2+\E(v)$. Define
\begin{equation}\label{eq:cone}
 \K=\{v\in\V:v\geq0\text{ a.e.}\}
\end{equation}
and the functional
\begin{equation}\label{eq:functional}
 J(v)=\frac12\E(v)+\int(a-F)v\dga,
 \qquad v\in\V.
\end{equation}
Then $J$ is finite and real-valued on $\V$, and $\K$ is nonempty,
convex, norm closed, and weakly closed in $\V$. Moreover:
\begin{enumerate}[label=\textup{(\roman*)},leftmargin=*]
 \item The infimum is finite and satisfies
 \[
  -\norm{F-m}_2^2\leq \iota:=\inf_{v\in\K}J(v)\leq0.
 \]
 \item Every fixed sublevel set
 $S_C:=\{v\in\K:J(v)\leq C\}$, $C\in\R$, is bounded in $\V$.
 More precisely, if $S_C$ is nonempty, then
 $M_C:=C+\norm{F-m}_2^2\geq0$ and
 \[
  \norm{v}_{\V}^2
  \leq\frac{M_C^2}{(a-m)^2}+8M_C,
  \qquad \forall v\in S_C.
 \]
 \item There is a unique $u\in\K$ such that $J(u)=\iota$.
\end{enumerate}
\end{lemma}

\begin{proof}
\textit{Step 1: Finiteness of the functional.}
Both terms in $J(v)$ are finite for every $v\in\V$. Indeed, by the
definition $\V=\Dom(L^{1/4})$, both $v$ and $L^{1/4}v$ belong to
$L^2(\ga)$, so
\[
 \E(v)=\ip{L^{1/4}v}{L^{1/4}v}
       =\norm{L^{1/4}v}_2^2<\infty.
\]
This does not require $v\in\U=\Dom(L^{1/2})$.
Since $F\in L^2(\ga)$ and $a1\in L^2(\ga)$,
Cauchy--Schwarz gives
\[
 \abs{\int(a-F)v\dga}
 \leq\int\abs{a-F}\abs{v}\dga
 \leq\norm{a1-F}_2\norm{v}_2<\infty.
\]
Since $\norm{v}_2\leq\norm{v}_{\V}$, this estimate also shows that the
linear term is a continuous linear functional on $\V$.
Thus $J$ is well defined and real-valued on all of $\V$, and in
particular on the admissible cone $\K$.

\medskip
\noindent\textit{Step 2: The admissible cone is weakly closed.}
The cone $\K$ contains $0$ and is convex.
To prove that $\K$ is weakly closed in $\V$, observe that for $v\in\V$,
\[
 v\in\K
 \quad\Longleftrightarrow\quad
 \int v\phi\dga\geq0
 \quad\text{for every }\phi\in L^2(\ga)\text{ with }\phi\geq0.
\]
The forward implication is obvious. For the converse, choose
$\phi=v^-=-\min \{v,0\}\in L^2(\ga)$. Then
\[
 0\leq\int v v^-\dga=-\norm{v^-}_2^2,
\]
so $v^-=0$ almost everywhere and $v\geq0$ almost everywhere.
Next, for each fixed $\phi\in L^2(\ga)$, define the map
$\Lambda_\phi:\V\to\R$ by
\[
 \Lambda_\phi(v):=\int v\phi\dga,\qquad v\in\V.
\]
This map is linear and continuous since 
\[
 |\Lambda_\phi(v)|
 \leq\norm{v}_2\norm{\phi}_2
 \leq\norm{v}_{\V}\norm{\phi}_2.
\]
For each nonnegative $\phi\in L^2(\ga)$, set
\[
 H_\phi:=\{v\in\V:\Lambda_\phi(v)\geq0\}
        =\Lambda_\phi^{-1}([0,\infty)).
\]
By definition of the weak topology, $\Lambda_\phi$ is weakly continuous.
Thus $H_\phi$ is weakly closed as the inverse image of the closed
interval $[0,\infty)$. The characterization of nonnegativity above gives
\[
 \K=\bigcap_{\substack{\phi\in L^2(\ga)\\\phi\geq0}}H_\phi.
\]
Hence $\K$ is weakly closed, being an intersection of weakly closed sets.

In particular, if $v_k\in\K$ and $v_k\rightharpoonup v$ in $\V$, then
for every nonnegative $\phi\in L^2(\ga)$,
\[
 \Lambda_\phi(v)
 =\lim_{k\to\infty}\Lambda_\phi(v_k)\geq0.
\]
Hence $v\in\K$, without requiring pointwise convergence. Since norm
convergence implies weak convergence, $\K$ is also norm closed.

\medskip
\noindent\textit{Step 3: A lower bound for $J$ and finiteness of its infimum.}
For an arbitrary $v\in\K$, put $c=\me(v)\geq0$ and $w=v-c1$.
Then $w\in\V$, $\me(w)=0$, and $v=c1+w$. Since $L^{1/4}1=0$,
we have $L^{1/4}w=L^{1/4}v$ and hence $\E(w)=\E(v)$.
The fractional spectral gap \eqref{eq:fractional-gap} therefore gives
$\norm{w}_2\leq\E(v)^{1/2}$. Since $m=\int F\dga$,
\[
 \begin{aligned}
 \int(a-F)v\dga
 &=\int(a-F)(c+w)\dga\\
 &=c\int(a-F)\dga+a\int w\dga-\int Fw\dga\\
 &=(a-m)c-\int(F-m)w\dga.
 \end{aligned}
\]
Here we used $\ga(\R^n)=1$, $\int w\dga=0$, and therefore
$\int Fw\dga=\int(F-m)w\dga$. By Cauchy--Schwarz and
\eqref{eq:fractional-gap},
\[
 \left|\int(F-m)w\dga\right|
 \leq\norm{F-m}_2\norm{w}_2
 \leq\norm{F-m}_2\E(v)^{1/2}.
\]
Substitution into the definition of $J$ gives the first inequality
below. For the second, Young's inequality gives
\[
 \norm{F-m}_2\E(v)^{1/2}
 \leq\norm{F-m}_2^2+\frac14\E(v).
\]
Consequently,
\begin{align}
 J(v)
 &\geq\frac12\E(v)-\norm{F-m}_2\E(v)^{1/2}+(a-m)c\notag\\
 &\geq\frac14\E(v)+(a-m)c-\norm{F-m}_2^2.
 \label{eq:coercivity}
\end{align}
Since $\E(v)\geq0$, $c\geq0$, and $a-m>0$, this inequality gives
$J(v)\geq-\norm{F-m}_2^2$ on $\K$.  Also $0\in\K$ and $J(0)=0$, so
the infimum is finite and satisfies
\[
 -\norm{F-m}_2^2\leq \iota:=\inf_{v\in\K}J(v)\leq0.
\]
\medskip
\noindent\textit{Step 4: An explicit bound for every fixed sublevel set.}
Fix $C\in\R$ and consider the set
\[
 S_C:=\{v\in\K:J(v)\leq C\}.
\]
We prove that all functions in $S_C$ have a common bound on their
$\V$-norm.  If $S_C$ is empty, there is nothing to prove.  Otherwise,
take an arbitrary $v\in S_C$ and put
$M_C=C+\norm{F-m}_2^2$.  With $c=\me(v)\geq0$,
\eqref{eq:coercivity} gives
\[
 0\leq\frac14\E(v)+(a-m)c\leq M_C.
\]
Since $a-m>0$, both terms are nonnegative.  Hence $M_C\geq0$ and
\[
 \E(v)\leq4M_C,
 \qquad
 0\leq c\leq\frac{M_C}{a-m}.
\]
Writing $v=c1+w$ as above, the mean-zero function $w$ is orthogonal
to $1$ in $L^2(\ga)$.  Pythagoras and \eqref{eq:fractional-gap} give
\[
 \norm{v}_2^2=c^2+\norm{w}_2^2
 \leq c^2+\E(v).
\]
Combining these estimates yields the explicit bound
\[
 \begin{aligned}
 \norm{v}_{\V}^2
 &=\norm{v}_2^2+\E(v)\leq c^2+2\E(v)\\
 &\leq\frac{M_C^2}{(a-m)^2}+8M_C.
 \end{aligned}
\]
The bound depends only on $C$ and the fixed data, not on $v$. Thus
every fixed sublevel set is bounded; no upper bound for $J$ on all
of $\K$ has been assumed.

\medskip
\noindent\textit{Step 5: Choosing a bounded minimizing sequence.}
By the definition of the finite infimum $\iota$, for each $k\geq1$ we may
choose $v_k\in\K$ such that
\[
 \iota\leq J(v_k)<\iota+\frac1k\leq\frac1k\leq1,
\]
where the last inequalities use $\iota\leq0$.  Hence $J(v_k)\to \iota$, so
$(v_k)$ is a minimizing sequence entirely contained in $S_1$.
The sublevel-set bound just proved, with $C=1$, gives
a common bound for $\norm{v_k}_{\V}$.  This construction does not
presuppose that the infimum is attained; it uses only the definition
of the infimum.

\medskip
\noindent\textit{Step 6: Weak convergence and attainment of the infimum.}
The space $\V$ is Hilbert, hence reflexive. The bounded sequence from
Step 5 therefore has a weakly convergent subsequence: for some $u\in\V$,
\[
 v_{k_j}\rightharpoonup u\quad\text{in }\V.
\]
Step 2 shows that $u\in\K$. To pass to the limit in the energy, note
that $L^{1/4}:\V\to L^2(\ga)$ is bounded and linear, because
\[
 \norm{L^{1/4}v}_2\leq\norm{v}_{\V}.
\]
Thus $L^{1/4}v_{k_j}\rightharpoonup L^{1/4}u$ in $L^2(\ga)$.
The weak lower semicontinuity of the Hilbert-space norm gives
\[
 \E(u)=\norm{L^{1/4}u}_2^2
 \leq\liminf_{j\to\infty}\norm{L^{1/4}v_{k_j}}_2^2
 =\liminf_{j\to\infty}\E(v_{k_j}).
\]
The linear term is weakly continuous: the Cauchy--Schwarz estimate in
Step 1 proves that $v\mapsto\int(a-F)v\dga$ is a continuous linear
functional on $\V$. Consequently,
\[
 \int(a-F)v_{k_j}\dga\longrightarrow\int(a-F)u\dga.
\]
Combining the two limits and using $J(v_{k_j})\to \iota$ yields
\[
 J(u)\leq\liminf_{j\to\infty}J(v_{k_j})=\iota.
\]
Since $u\in\K$, the definition of $\iota$ also gives $\iota\leq J(u)$.
Therefore $J(u)=\iota$, and the infimum is attained. This argument uses
weak compactness of bounded sequences in $\V$, not a compact embedding
into another function space; no compactness uniform in $n$ is needed.

\medskip
\noindent\textit{Step 7: Uniqueness, including the constant part.}
Suppose that $u_1,u_2\in\K$ are minimizers. Their midpoint belongs to
$\K$ by convexity. The parallelogram identity for the quadratic energy
and linearity of the remaining term give
\[
 J\left(\frac{u_1+u_2}{2}\right)
 =\frac{J(u_1)+J(u_2)}2-\frac18\E(u_1-u_2).
\]
Since $J(u_1)=J(u_2)=\iota$ and their midpoint is admissible,
minimality gives
\[
 \iota\leq J\left(\frac{u_1+u_2}{2}\right)
 =\iota-\frac18\E(u_1-u_2).
\]
Hence $\E(u_1-u_2)\leq0$. On the other hand,
\[
 \E(u_1-u_2)=\norm{L^{1/4}(u_1-u_2)}_2^2\geq0,
\]
so equality must hold. By \eqref{eq:form-series}, every nonconstant
Hermite coefficient of $u_1-u_2$ is zero, so $u_1-u_2=d1$ for some
$d\in\R$. Adding a constant does not change the energy, and hence
\[
 0=J(u_1)-J(u_2)
   =\int(a-F)(u_1-u_2)\dga=(a-m)d.
\]
Since $a-m>0$, we obtain $d=0$ and $u_1=u_2$ almost everywhere.
Thus the energy determines the mean-zero part of the minimizer, while
the strict height inequality determines its constant part.
\end{proof}

The minimizer is now constructed in $\V$. The next two lemmas
control its multiplier; \cref{subsec:upgrade} then supplies the
stronger operator-domain regularity.

\subsection{The variational inequality and the multiplier}
\label{subsec:duality}

\begin{lemma}[Variational inequality and obstacle multiplier]
\label{lem:obstacle-multiplier}
Let $F$, $a$, $m$, $\V$, $\K$, and $J$ be as in
\cref{lem:obstacle-minimizer}, and let $u\in\K$ be the unique
minimizer. Denote by $\V^*$ the continuous dual of $\V$, with
duality pairing $\dual{\Lambda}{\phi}=\Lambda(\phi)$.
Then $u$ satisfies
\begin{equation}\label{eq:VI}
 \E(u,v-u)\geq\int(F-a)(v-u)\dga,
 \qquad v\in\K.
\end{equation}
The formula
\begin{equation}\label{eq:eta-functional}
 \dual{\eta}{\phi}
    =\E(u,\phi)-\int(F-a)\phi\dga,
 \qquad \phi\in\V,
\end{equation}
defines the obstacle multiplier $\eta\in\V^*$, with
\[
 \norm{\eta}_{\V^*}\leq\E(u)^{1/2}+\norm{F-a}_2.
\]
The multiplier is the Fr\'echet derivative $J'(u)$ of the functional
$J:\V\to\R$ at $u$. This derivative is an element of $\V^*$;
its value on a direction $\phi\in\V$ is
\[
 J'(u)[\phi]
 =\left.\frac{d}{dt}J(u+t\phi)\right|_{t=0}
 =\dual{\eta}{\phi}.
\]
Thus $J'(u)=\eta$ is an equality of continuous linear functionals,
not a pointwise derivative. The multiplier is positive on nonnegative tests:
\begin{equation}\label{eq:eta-positive}
 \dual{\eta}{\phi}\geq0,
 \qquad \phi\in\V,\quad\phi\geq0.
\end{equation}
Moreover,
\[
 \dual{\eta}{u}=0,
 \qquad \dual{\eta}{1}=a-m>0.
\]
The variational inequality is equivalent to positivity
of $\eta$ together with $\dual{\eta}{u}=0$. At this stage $\eta$ is
a functional; no pointwise values or $L^2$ representation are assumed.
\end{lemma}

\begin{proof}
\noindent\textit{Step 1: The duality pairing and the $L^2$ embedding.}
An element $\Lambda\in\V^*$ is a linear map $\V\to\R$ whose
dual norm
\[
 \norm{\Lambda}_{\V^*}
   =\sup_{\norm{\phi}_{\V}\leq1}|\Lambda(\phi)|
\]
is finite. The pairing $\dual{\Lambda}{\phi}$ means evaluation of
this functional at $\phi$, not necessarily an integral.
Every $h\in L^2(\ga)$ induces a functional $\ell_h\in\V^*$ by
\[
 \dual{\ell_h}{\phi}=\ip{h}{\phi}=\int h\phi\dga,
 \qquad \phi\in\V,
\]
since Cauchy--Schwarz and $\norm{\phi}_2\leq\norm{\phi}_{\V}$ give
$|\dual{\ell_h}{\phi}|\leq\norm{h}_2\norm{\phi}_{\V}$.
This embedding is injective because $\V$ is dense in $L^2(\ga)$,
but not every element of $\V^*$ has such an $L^2$ density.

\medskip
\noindent\textit{Step 2: The form operator and the definition of $\eta$.}
We know $u\in\V=\Dom(L^{1/4})$, but not yet
$u\in\U=\Dom(B)$. Thus $Bu$ cannot yet be treated as an $L^2$
function. Instead, define $\mathcal B:\V\to\V^*$ by
\[
 \dual{\mathcal B v}{\phi}=\E(v,\phi),
 \qquad v,\phi\in\V.
\]
Cauchy--Schwarz gives
\[
 |\E(v,\phi)|
 \leq\E(v)^{1/2}\E(\phi)^{1/2}
 \leq\E(v)^{1/2}\norm{\phi}_{\V},
\]
so $\mathcal B$ is bounded, with
$\norm{\mathcal B v}_{\V^*}\leq\E(v)^{1/2}\leq\norm{v}_{\V}$.
For $v\in\U$, the form identity in
\cref{rem:quarter-power-domain} gives $\mathcal B v=\ell_{Bv}$.
Since $F-a1\in L^2(\ga)$, Step 1 also gives
$\ell_{F-a}\in\V^*$. Hence \eqref{eq:eta-functional} defines
\[
 \eta=\mathcal B u-\ell_{F-a}\in\V^*,
 \qquad
 \norm{\eta}_{\V^*}\leq\E(u)^{1/2}+\norm{F-a}_2.
\]
The shorthand $\eta=Bu+a-F$ refers to this equality of functionals,
not yet to a pointwise or $L^2$ equality.

\medskip
\noindent\textit{Step 3: First variation and the variational inequality.}
Bilinearity of $\E$ gives, for every $h\in\V$,
\[
 J(u+h)-J(u)=\dual{\eta}{h}+\frac12\E(h).
\]
The Fr\'echet derivative is the continuous linear functional that
approximates this change with error negligible relative to
$\norm{h}_{\V}$. Here $\E(h)\leq\norm{h}_{\V}^2$ gives, for $h\ne0$,
\[
 \frac{\abs{J(u+h)-J(u)-\dual{\eta}{h}}}{\norm{h}_{\V}}
 =\frac{\E(h)}{2\norm{h}_{\V}}
 \leq\frac12\norm{h}_{\V}\longrightarrow0
 \quad\text{as }h\to0\text{ in }\V.
\]
This proves $J'(u)=\eta$ in $\V^*$; taking $h=t\phi$ also gives
the directional derivative formula in the statement.
For $v\in\K$ and $0<\theta\leq1$, convexity gives
$u+\theta(v-u)\in\K$. Minimality and the same expansion imply
\[
 0\leq\frac{J(u+\theta(v-u))-J(u)}{\theta}
   =\dual{\eta}{v-u}+\frac\theta2\E(v-u).
\]
Letting $\theta\downarrow0$ yields
\[
 \dual{\eta}{v-u}\geq0,\qquad v\in\K.
\]
By \eqref{eq:eta-functional}, this is exactly \eqref{eq:VI}.
The derivative need not vanish in every direction, because $u$
minimizes $J$ only over the constrained set $\K$.

\medskip
\noindent\textit{Step 4: Positivity and dual complementarity.}
For a nonnegative $\phi\in\V$, the competitor $v=u+\phi$ belongs
to $\K$, so Step 3 gives $\dual{\eta}{\phi}\geq0$.
This proves \eqref{eq:eta-positive}; the notation $\eta\geq0$ in
$\V^*$ means positivity on nonnegative tests, not a pointwise claim.
The competitors $v=0$ and $v=2u$ give, respectively,
\[
 -\dual{\eta}{u}\geq0,
 \qquad \dual{\eta}{u}\geq0,
\]
and hence $\dual{\eta}{u}=0$. Conversely, positivity and this equality
give $\dual{\eta}{v-u}=\dual{\eta}{v}\geq0$ for every $v\in\K$,
proving the stated equivalence. Finally, $1\in\V$ and $\E(u,1)=0$,
so \eqref{eq:eta-functional} gives
\[
 \dual{\eta}{1}=\E(u,1)-\int(F-a)\dga=a-m>0.
\]

\end{proof}

The Riesz representation theorem on $\V$ alone would give a vector
$z\in\V$ with $\dual{\eta}{\phi}=\ip{z}{\phi}+\E(z,\phi)$, not
an $L^2$ density for $\eta$. To obtain such a density we need the
stronger estimate $|\dual{\eta}{\phi}|\leq C\norm{\phi}_2$.
The next lemma supplies exactly the domination from which this
stronger bound follows.

\subsection{The Lewy--Stampacchia estimate on the full form domain}
\label{subsec:LS}

We first explain the statement and its role. In an obstacle problem,
the minimizer is constrained to remain above an obstacle; here the
obstacle is zero and the constraint is $u\geq0$. The multiplier
$\eta$ measures the reaction exerted by this constraint. Formally,
\eqref{eq:eta-functional} reads
\[
 \eta=Bu+a-F.
\]
The classical Lewy--Stampacchia estimate appears in \cite{LS70}. For example, if $\Omega$
is a bounded open set, $\psi\in C_c^\infty(\Omega)$, and $w$ minimizes
$\tfrac12\int_\Omega|\nabla v|^2\,dx$ over
$v\in W^{1,2}_0(\Omega)$ with $v\geq\psi$, the estimate reads
\[
 0\leq-\Delta w\leq(-\Delta\psi)^+.
\]
It bounds the reaction $-\Delta w$ from above by the positive part
of the operator applied to the obstacle. Abstract and fractional
versions are given by Gigli--Mosconi
\cite[Theorem~2.4 and Section~3.2]{GM15}.

Here positivity of $\eta$ already follows from \eqref{eq:eta-positive};
the additional Lewy--Stampacchia conclusion is the upper bound by
$(a-F)^+$. Its precise dual-order formulation is proved directly in
\cref{lem:LS-order} below. We include the proof rather than invoke an abstract obstacle theorem.
After $\eta$ is represented by a function in the next subsection, this
becomes the pointwise estimate
\[
 0\leq\eta\leq(a-F)^+.
\]
Equivalently, because $Bu=F-a+\eta$, it gives
\[
 F-a\leq Bu\leq(F-a)^+.
\]
This estimate is decisive for two reasons. First, the domination
$\eta\leq(a-F)^+\leq a$ makes the multiplier an $L^2$ function rather
than merely an element of $\V^*$. Second, after setting $g=a-\eta$, it
immediately yields the desired bounds
$\min\{F,a\}\leq g\leq a$. The upper estimate is proved by perturbing
$u$ downward in a nonnegative direction and truncating the perturbation
at the obstacle; this is the purpose of the competitor below.

\begin{lemma}[Lewy--Stampacchia order estimate]
\label{lem:LS-order}
Let $u\in\K$ be the minimizer of $J$, let $\eta\in\V^*$ be defined by
\eqref{eq:eta-functional}, and set $\kappa_F=(a-F)^+$. Then
\begin{equation}\label{eq:eta-order-V}
 0\leq\dual{\eta}{\phi}\leq\int \kappa_F\phi\dga,
 \qquad \phi\in\V,\quad\phi\geq0.
\end{equation}
Thus the multiplier $\eta$ is bounded above, in the order of $\V^*$, by the
functional induced by the function $(a-F)^+$.
\end{lemma}

\begin{proof}
Fix $\phi\in\V$, $\phi\geq0$, and $\varepsilon>0$. Set
\[
 z=u-\varepsilon\phi,
 \qquad v=z^+\in\K.
\]
Indeed, $z\in\V$, and \cref{lem:lattice} shows that $z^+\in\V$;
since $z^+\geq0$, it follows that $v=z^+\in\K$. Thus $v$ is an
admissible competitor in \eqref{eq:VI}. Using
$z=z^+-z^-$ and $u=z+\varepsilon\phi$, we obtain
\[
 v-u=z^+-(z+\varepsilon\phi)=z^--\varepsilon\phi.
\]
Substitution into \eqref{eq:VI} yields
\[
 \E(u,z^-)-\varepsilon\E(u,\phi)
 \geq\int(F-a)z^-\dga
       -\varepsilon\int(F-a)\phi\dga.
\]
Moving the terms containing $\varepsilon$ to the left and using
\eqref{eq:eta-functional}, we find
\begin{equation}\label{eq:LS-start}
 \varepsilon\dual{\eta}{\phi}
 \leq\E(u,z^-)+\int(a-F)z^-\dga.
\end{equation}

We first estimate the energy term on the right. Since
$u=z+\varepsilon\phi$, bilinearity gives
$\E(u,z^-)=\E(z,z^-)+\varepsilon\E(\phi,z^-)$. The sign estimate
\eqref{eq:negative-part-sign} controls the first term, while
Cauchy--Schwarz for the form $\E$ controls the second:
\begin{align}
 \E(u,z^-)
 &=\E(z,z^-)+\varepsilon\E(\phi,z^-)\notag\\
 &\leq-\E(z^-)
       +\varepsilon\E(\phi)^{1/2}\E(z^-)^{1/2}\notag\\
 &\leq\frac{\varepsilon^2}{4}\E(\phi).
 \label{eq:LS-energy}
\end{align}
For the last step, set
$p=\E(\phi)^{1/2}$ and $q=\E(z^-)^{1/2}$. Completing the square gives
\[
 -q^2+\varepsilon pq
 =-\left(q-\frac{\varepsilon p}{2}\right)^2
   +\frac{\varepsilon^2p^2}{4}
 \leq\frac{\varepsilon^2}{4}\E(\phi),
\]
which proves \eqref{eq:LS-energy}.

We next estimate the load term. Because $u\geq0$ and $\phi\geq0$,
\[
 0\leq z^-=(\varepsilon\phi-u)^+
 \leq\varepsilon\phi
 \quad\text{a.e.}
\]
Recall that $\kappa_F=(a-F)^+$. We have $a-F\leq \kappa_F$ and $\kappa_F\geq0$. Therefore
\[
 (a-F)z^-\leq kz^-\leq\varepsilon \kappa_F\phi
 \quad\text{a.e.},
\]
and integration gives
\begin{equation}\label{eq:LS-load}
 \int(a-F)z^-\dga\leq\varepsilon\int \kappa_F\phi\dga.
\end{equation}
The integral on the right is finite: since $F\geq0$, we have
$0\leq \kappa_F=(a-F)^+\leq a$, and $\phi\in L^2(\ga)\subseteq L^1(\ga)$.

Combining \eqref{eq:LS-start}, \eqref{eq:LS-energy}, and
\eqref{eq:LS-load}, we obtain
\[
 \varepsilon\dual{\eta}{\phi}
 \leq\frac{\varepsilon^2}{4}\E(\phi)
      +\varepsilon\int \kappa_F\phi\dga.
\]
Since $\varepsilon>0$, division by $\varepsilon$ gives
\[
 \dual{\eta}{\phi}
 \leq\frac{\varepsilon}{4}\E(\phi)+\int \kappa_F\phi\dga.
\]
The test function $\phi$ is fixed, so $\E(\phi)<\infty$; letting
$\varepsilon\downarrow0$ removes the first term. Combining the resulting
upper bound with the positivity statement \eqref{eq:eta-positive} gives
exactly \eqref{eq:eta-order-V}. The estimate has therefore been obtained
for \emph{every} nonnegative form-domain test, not merely for smooth
compactly supported functions.
\end{proof}

\subsection{Representing the multiplier and upgrading the domain}
\label{subsec:upgrade}

So far $u\in\V=\Dom(L^{1/4})$ and $\eta\in\V^*$: neither $Bu$
nor $\eta$ has yet been identified as an $L^2$ function. This subsection
makes two upgrades. The order estimate first gives an $L^2$ density
for $\eta$. The resulting form equation then puts $u$ in the smaller
operator domain $\Dom(B)=\U$, where $Bu$ and $Du$ are well defined
in $L^2$.

\begin{lemma}[$L^2$ representation and operator-domain upgrade]
\label{lem:multiplier-domain-upgrade}
Let $u$ and $\eta$ be as in \cref{lem:obstacle-multiplier}, and set
$\kappa_F=(a-F)^+$. By \cref{lem:LS-order}, the dual-order estimate
\eqref{eq:eta-order-V} holds.
Then $\eta$ has a unique representing function in $L^2(\ga)$, also
denoted by $\eta$, satisfying
\[
 \dual{\eta}{\phi}=\int\eta\phi\dga\quad(\phi\in\V),
 \qquad 0\leq\eta\leq \kappa_F\leq a\quad\text{a.e.}
\]
The function $g=a-\eta$ belongs to $L^2(\ga)$ and satisfies
$\min\{F,a\}\leq g\leq a$ almost everywhere. Moreover,
\[
 u\in\Dom(B)=\U=W^{1,2}(\ga),
 \qquad Bu=F-g\in L^2(\ga),
\]
and $F-a\leq Bu\leq(F-a)^+$ almost everywhere.
\end{lemma}

\begin{proof}
\noindent\textit{Step 1: Obtain an $L^2$ density for the functional $\eta$.}
For any real $\psi\in\V$, the lattice property gives
$\psi^+,\psi^-,|\psi|\in\V$. Since
$\psi=\psi^+-\psi^-$ and $\eta$ is positive on nonnegative tests,
$|\dual{\eta}{\psi}|\leq\dual{\eta}{\psi^+}+\dual{\eta}{\psi^-}
=\dual{\eta}{|\psi|}$. The order estimate and $\kappa_F\leq a$ therefore give
\begin{equation}\label{eq:eta-L2-bound}
 |\dual{\eta}{\psi}|
 \leq\dual{\eta}{|\psi|}
 \leq\int \kappa_F|\psi|\dga
 \leq a\norm{\psi}_1\leq a\norm{\psi}_2.
\end{equation}
The last inequality uses $\ga(\R^n)=1$. Crucially, this bound uses
$\norm{\psi}_2$, not $\norm{\psi}_{\V}$.
Finite Hermite sums show that $\V$ is dense in $L^2$. Consequently,
$\eta$ extends uniquely to a bounded functional on $L^2$: for an
$L^2$ approximation $\psi_j\in\V$ to $f$, define its value at $f$
as $\lim_j\dual{\eta}{\psi_j}$. The bound above ensures that this
limit exists and is independent of the approximation. The $L^2$
Riesz theorem then supplies the unique density, which we denote
again by $\eta$.

\medskip
\noindent\textit{Step 2: Convert dual order into pointwise bounds.}
An indicator $\one_A$ need not belong to $\V$, so it cannot be used
directly in \eqref{eq:eta-order-V}. Instead, for any measurable
$A\subseteq\R^n$, use $P_s\one_A$ with $s>0$. The positivity-preserving
property established in \cref{subsec:Markov} makes this test nonnegative.
The spectral action $P_sh_\alpha=e^{-s|\alpha|}h_\alpha$ recorded there,
together with the domain formula \eqref{eq:spectral-domains}, gives
\[
 \norm{L^{1/4}P_s\one_A}_2^2
 \leq\left(\sup_{\ell\geq1}\ell^{1/2}e^{-2s\ell}\right)
      \norm{\one_A}_2^2<\infty.
\]
Thus $P_s\one_A\in\V$. Strong continuity gives
$P_s\one_A\to\one_A$ in $L^2$ as $s\downarrow0$. Apply
\eqref{eq:eta-order-V} to these tests and pass to the limit using
$\kappa_F,\eta\in L^2$:
\[
 0\leq\int_A\eta\dga\leq\int_A \kappa_F\dga.
\]
Taking $A=\{\eta<0\}$ and $A=\{\eta>\kappa_F\}$ rules out either set
having positive measure. Hence
\begin{equation}\label{eq:eta-pointwise}
 0\leq\eta\leq(a-F)^+\leq a\quad\text{a.e.}
\end{equation}

\medskip
\noindent\textit{Step 3: Construct $g$ and obtain an $L^2$ form equation.}
Define
\begin{equation}\label{eq:g-definition}
 g=a-\eta,
 \qquad b=F-g=F-a+\eta\in L^2(\ga).
\end{equation}
Because $a1$, $F$, and the density $\eta$ all lie in $L^2$, both
$g$ and $b$ are genuine $L^2$ functions. Substituting the density into
\eqref{eq:eta-functional} yields
\begin{equation}\label{eq:form-equation}
 \E(u,\phi)=\ip{b}{\phi},\qquad\phi\in\V.
\end{equation}

\medskip
\noindent\textit{Step 4: Upgrade $u$ to the operator domain.}
Every Hermite basis function belongs to $\V$. Testing
\eqref{eq:form-equation} against $h_0=1$ gives $\ip{b}{1}=0$;
for $\alpha\ne0$, testing against $h_\alpha$ gives
\[
 \ip{b}{h_\alpha}=\E(u,h_\alpha)
        =|\alpha|^{1/2}u_\alpha.
\]
Parseval's identity therefore gives
\begin{equation}\label{eq:upgrade-series}
 \sum_{\alpha\ne0}|\alpha|\,|u_\alpha|^2
     =\sum_{\alpha\ne0}|\ip{b}{h_\alpha}|^2
     =\norm{b}_2^2<\infty.
\end{equation}
This is exactly the domain condition for $u\in\Dom(L^{1/2})$ in
\eqref{eq:spectral-domains}. The Hermite coefficients of $Bu$ and
$b$ agree, including the zero coefficient, so $Bu=b=F-g$ in $L^2$.
By \eqref{eq:sqrt-domain}, $u\in W^{1,2}(\ga)$ and $Du\in L^2$.
The identity $\Rg Bu=Du$ in \eqref{eq:R-B} is now justified.

\medskip
\noindent\textit{Step 5: Recover the pointwise operator and good-part bounds.}
Since $Bu=F-a+\eta$ and $0\leq\eta\leq(a-F)^+$,
\begin{equation}\label{eq:LS-B}
 F-a\leq Bu\leq F-a+(a-F)^+=(F-a)^+.
\end{equation}
This is the full Lewy--Stampacchia estimate as a pointwise inequality.
Likewise, $g=a-\eta$ gives
\[
 \min\{F,a\}=a-(a-F)^+\leq g\leq a.
\]
In particular, $0\leq g\leq a$, as required for the good part of
\cref{thm:decomposition}.
\end{proof}

\subsection{Conservation of Gaussian mass}
\label{subsec:mass}

The constant function is an admissible form-domain test. Equation
\eqref{eq:form-equation} with $\phi=1$ yields
\begin{equation}\label{eq:mass-proof}
 0=\E(u,1)=\int(F-g)\dga.
\end{equation}
Therefore $\int g\dga=m$. All terms are integrable by
\eqref{eq:L2intoL1}. There is no limiting cutoff at spatial infinity
and no estimate for the mass of a family of balls.
Equivalently, $\int\eta\dga=a-m$.

\subsection{Complementarity and the exceptional set}
\label{subsec:complementarity}
For a measurable representative of $u$, set $\Omega=\{u>0\}$.
We now prove that $g=a$ almost everywhere on $\Omega$ and that
$a\ga(\Omega)\leq m$. These are precisely
\eqref{eq:g-contact} and \eqref{eq:Omega-bound} in
\cref{thm:decomposition}. The key is \emph{complementarity}:
$\eta u=0$ almost everywhere, so the potential $u$ and the obstacle
reaction $\eta$ cannot both be positive at the same point, except
on a null set.

By \cref{lem:obstacle-multiplier}, $\dual{\eta}{u}=0$. Now that
$\eta$ is an $L^2$ function, this is an integral identity. In particular,
\eqref{eq:eta-functional} gives
\begin{equation}\label{eq:energy-equality}
 \E(u)=\int(F-a)u\dga,
\end{equation}
and $a-g=\eta$ yields
\begin{equation}\label{eq:complementarity}
 0=\int\eta u\dga=\int(a-g)u\dga.
\end{equation}
The product $\eta u$ is integrable by Cauchy--Schwarz, since
$\eta,u\in L^2(\ga)$. It is also nonnegative. A nonnegative
integrable function with zero integral vanishes almost everywhere,
so $\eta u=0$ a.e. On $\Omega$, where $u>0$, this forces $\eta=0$;
hence $g=a-\eta=a$ almost everywhere on $\Omega$, proving
\eqref{eq:g-contact}.

\medskip
\noindent\textit{The exceptional-set bound.}
Because $g=a$ on $\Omega$, $g\geq0$ everywhere up to a null set,
and $\int g\dga=m$ by \eqref{eq:mass-proof},
\begin{equation}\label{eq:exceptional-proof}
 a\ga(\Omega)=\int_\Omega g\dga
      \leq\int g\dga=m.
\end{equation}
This proves \eqref{eq:Omega-bound}. Only measurability of $\Omega$
is needed, not openness or boundary regularity. The remaining
gradient-support assertion \eqref{eq:gradient-support} of
\cref{thm:decomposition} is proved in \cref{subsec:locality}.

\paragraph{Uniqueness of the pair.}
Suppose $\widetilde u\in\Dom(B)$ and $\widetilde g\in L^2$ satisfy
$\widetilde u\geq0$, $F=\widetilde g+B\widetilde u$,
$\widetilde g\leq a$, and $\widetilde g=a$ on $\{\widetilde u>0\}$.
Put $\widetilde\eta=a-\widetilde g\geq0$. Then
$\widetilde\eta\widetilde u=0$, and for every $v\in\K$,
\[
 \E(\widetilde u,v-\widetilde u)
       -\int(F-a)(v-\widetilde u)\dga
 =\int\widetilde\eta v\dga\geq0.
\]
The quadratic expansion of $J$ shows that $\widetilde u$ minimizes $J$.
Uniqueness in \cref{lem:obstacle-minimizer} gives $\widetilde u=u$,
and then $\widetilde g=F-Bu=g$. This proves the uniqueness asserted in
\cref{thm:decomposition}.

\subsection{Why the gradient vanishes on the zero set}
\label{subsec:locality}

Choose a nonnegative measurable representative of $u$ and set
$\Omega=\{u>0\}$. Its complement is the contact set with the zero
obstacle. We prove the last assertion
\eqref{eq:gradient-support} of \cref{thm:decomposition}:
\[
 Du=0\quad\text{a.e. on }\Omega^c=\{u=0\}.
\]
Here $u\geq0$ and $u\in\U=W^{1,2}(\ga)$ by
\cref{lem:multiplier-domain-upgrade}, and $D=2^{-1/2}\nabla$ is the
normalized ordinary weak gradient. The connection to the Gaussian Riesz
transform is the identity already proved in \eqref{eq:R-B}:
\[
 \Rg(F-g)=\Rg Bu=Du=\frac1{\sqrt2}\nabla u.
\]
The new conclusion is that this vector field vanishes outside $\Omega$.
Thus the argument concerns locality of the ordinary gradient, not of
the nonlocal operator $B=L^{1/2}$.

\medskip
\noindent\textit{Step 1: Apply the chain rule to smooth contractions of $u$.}
For $\delta>0$, define
\[
 q_\delta(t)=\sqrt{t^2+\delta^2}-\delta,
 \qquad q_\delta'(t)=\frac{t}{\sqrt{t^2+\delta^2}}.
\]
The function $q_\delta$ is smooth, vanishes at zero, and is
$1$-Lipschitz because $|q_\delta'|\leq1$. Hence
$q_\delta(u)\in W^{1,2}(\ga)$, and the Sobolev chain rule gives
\begin{equation}\label{eq:chain-rule}
 D(q_\delta(u))
   =\frac{u}{\sqrt{u^2+\delta^2}}\,Du.
\end{equation}
To justify the chain rule, approximate $u$ in $W^{1,2}(\ga)$ by smooth
functions $u_j$, choosing a subsequence converging almost everywhere.
The Lipschitz bound gives $q_\delta(u_j)\to q_\delta(u)$ in $L^2$.
For the gradients, boundedness of $q_\delta'$ controls the term
$q_\delta'(u_j)(Du_j-Du)$, while continuity of $q_\delta'$ and
dominated convergence control $(q_\delta'(u_j)-q_\delta'(u))Du$.
Both terms tend to zero in $L^2$. Closedness of $D$, established in
\cref{sec:spectral}, then proves \eqref{eq:chain-rule}.

\medskip
\noindent\textit{Step 2: Identify the limits of the functions and gradients.}
Since $u\geq0$,
\[
 0\leq q_\delta(u)\leq u,
 \qquad q_\delta(u)\longrightarrow u
 \quad\text{a.e. as }\delta\downarrow0.
\]
As $u\in L^2$, dominated convergence yields $q_\delta(u)\to u$
in $L^2(\ga)$. The multiplier in \eqref{eq:chain-rule} satisfies
\[
 0\leq\frac{u}{\sqrt{u^2+\delta^2}}\leq1,
 \qquad
 \frac{u}{\sqrt{u^2+\delta^2}}\longrightarrow\one_\Omega
 \quad\text{a.e.}
\]
Indeed, its limit is $1$ where $u>0$ and it is identically $0$ where
$u=0$. Since $Du\in L^2(\ga;\R^n)$, dominated convergence gives
$D(q_\delta(u))\to\one_\Omega Du$ in $L^2(\ga;\R^n)$.

\medskip
\noindent\textit{Step 3: Use closedness of the gradient and localize the transform.}
The gradient $D$ is closed: if $v_j\in W^{1,2}(\ga)$,
$v_j\to v$ in $L^2$, and
$Dv_j\to G$ in $L^2(\ga;\R^n)$, then $v\in W^{1,2}(\ga)$ and
$Dv=G$. Applying this to $v_j=q_{\delta_j}(u)$ for any
$\delta_j\downarrow0$ and using Step 2 yields
\begin{equation}\label{eq:locality-proof}
 Du=\one_\Omega Du.
\end{equation}
Thus $Du=0$ almost everywhere on $\Omega^c$, proving
\eqref{eq:gradient-support} without any openness or boundary-regularity
assumption on $\Omega$. Combining this with \eqref{eq:R-B} and
$F=g+Bu$ gives
\[
 \Rg F=\Rg g+\frac1{\sqrt2}\nabla u,
 \qquad \Rg F=\Rg g\quad\text{a.e. on }\Omega^c.
\]
This is the exact localization needed for the endpoint argument:
the transformed bad part vanishes outside $\Omega$, even though
the bad function $Bu=F-g$ itself need not vanish there.
Together with the preceding subsections, this completes the proof of
\cref{thm:decomposition}.

\begin{remark}[What the proof does and does not use]
The only role of the spectral gap is to obtain a minimizer. No numerical
bound for its size is used afterward. The order estimates, mass identity,
and exceptional-set estimate have exact constants. The proof never
estimates a covering multiplicity, a Gaussian doubling ratio, or the
absolute size of the Mehler kernel.
\end{remark}

\section{The dimension-free endpoint argument}
\label{sec:endpoint}

\begin{theorem}[Dimension-free weak-type $(1,1)$ estimate]
\label{thm:endpoint}
The transform in \cref{def:riesz}, initially defined on real
$L^2(\ga)$, has a unique continuous linear extension from $L^1(\ga)$ to
$L^{1,\infty}(\ga;\R^n)$ satisfying
\begin{equation}\label{eq:endpoint-final}
 \sup_{t>0}t\,\ga\bigl(\{|\Rg f|_{\R^n}>t\}\bigr)
       \leq2\norm{f}_1.
\end{equation}
This estimate holds for every $n\geq1$.
\end{theorem}

We first prove the estimate for real $L^2$ inputs. The signed potential
and good function will be obtained by applying \cref{thm:decomposition}
to $f^+$ and $f^-$ at the same height, and subtracting the resulting
pairs. The $L^1$ extension is constructed afterwards.

\subsection{Real $L^2$ inputs: the signed decomposition}
\label{subsec:signed}

Let $f\in L^2(\ga;\R)$ and set $M=\norm{f}_1$. If $M=0$, the statement
is immediate. By \eqref{eq:small-level}, thresholds $0<t\leq M$ are
already controlled with constant $1$. Fix $t>M$.

Apply \cref{thm:decomposition} at height $a=t$ first with $F=f^+$
and then with $F=f^-$. This is the first construction step. It produces
nonnegative potentials $u_+,u_-\in\Dom(B)$ and good functions
$g_+,g_-\in L^2(\ga)\cap L^\infty(\ga)$ satisfying
\begin{equation}\label{eq:signed-pieces}
 f^+=g_++Bu_+,\qquad f^-=g_-+Bu_-.
\end{equation}
The notation $u_+$ and $u_-$ denotes two independently constructed
nonnegative potentials; it does \emph{not} assert that they are the
positive and negative parts of one function. The same caution applies
to $g_+$ and $g_-$.

We now perform the second construction step and define the signed
potential, signed good function, and exceptional set by
\begin{equation}\label{eq:signed-definition}
 \begin{split}
 g&=g_+-g_-,\qquad u=u_+-u_-,\\
 \Omega&=\{u_+>0\}\cup\{u_->0\}.
 \end{split}
\end{equation}
Because $\Dom(B)$ is a vector space, $u\in\Dom(B)$; similarly,
$g\in L^2(\ga)\cap L^\infty(\ga)$. Subtracting the two identities in
\eqref{eq:signed-pieces} and using linearity of $B$ gives
\[
 f=f^+-f^-=g_+-g_-+B(u_+-u_-)=g+Bu
 \quad\text{in }L^2(\ga).
\]
Thus \eqref{eq:signed-definition} constructs exactly the $u$ and $g$
needed in the remainder of \cref{sec:endpoint}. The exceptional-set estimates imply
\begin{equation}\label{eq:signed-Omega}
 \ga(\Omega)\leq\frac{\int f^+\dga+\int f^-\dga}{t}
                 =\frac{M}{t}.
\end{equation}
Moreover, $Du=0$ a.e. on $\Omega^c$.

Because $0\leq g_+,g_-\leq t$, their difference satisfies
\begin{equation}\label{eq:signed-sup}
 |g|\leq t,
\end{equation}
not merely $2t$. Mass preservation gives
\begin{equation}\label{eq:signed-L1}
 \norm{g}_1\leq\int(g_++g_-)\dga
       =\int(f^++f^-)\dga=M.
\end{equation}
The first step is an inequality: $g_+$ and $g_-$ need not have disjoint
supports. Combining the last two bounds,
\begin{equation}\label{eq:good-L2}
 \norm{g}_2^2\leq t\norm{g}_1\leq tM.
\end{equation}

\subsection{Exact localization of the transformed bad part}

All terms of $f=g+Bu$ are in $L^2$, and
$u\in\Dom(B)$. We may therefore apply the bounded $L^2$ operator
$\Rg$ and use \cref{prop:identities} to obtain
\begin{equation}\label{eq:transformed-decomposition}
 \Rg f=\Rg g+Du.
\end{equation}
Outside $\Omega$, the second term is zero a.e. Consequently,
\begin{equation}\label{eq:level-inclusion}
 \{|\Rg f|>t\}\subseteq
          \Omega\cup\{|\Rg g|>t\}
 \quad\text{up to a null set}.
\end{equation}

By \eqref{eq:signed-Omega}, Chebyshev's inequality,
\eqref{eq:R-L2}, and \eqref{eq:good-L2},
\begin{align}
 \ga(\{|\Rg f|>t\})
 &\leq\ga(\Omega)+t^{-2}\norm{\Rg g}_2^2\notag\\
 &\leq\frac{M}{t}+t^{-2}\norm{\Pi_0g}_2^2\notag\\
 &\leq\frac{M}{t}+t^{-2}\norm{g}_2^2\notag\\
 &\leq\frac{2M}{t}.
 \label{eq:main-bound}
\end{align}
Together with the small-level estimate, this proves
\eqref{eq:endpoint-final} for real $L^2$ functions.

\subsection{Construction of the $L^1$ extension}
\label{subsec:L1extension}

Let $f\in L^1(\ga;\R)$ and choose real $f_k\in L^2(\ga)$ with
$f_k\to f$ in $L^1$. For example, the pointwise truncations
$f_k=\max\{-k,\min\{f,k\}\}$ are bounded and hence in $L^2$.
For every $\varepsilon>0$, the estimate already proved implies
\begin{equation}\label{eq:Cauchy-measure}
 \ga\bigl(\{|\Rg f_k-\Rg f_\ell|>\varepsilon\}\bigr)
       \leq\frac{2}{\varepsilon}\norm{f_k-f_\ell}_1.
\end{equation}
Thus $\Rg f_k$ is Cauchy in measure.

Here is an explicit way to obtain its limit. Choose a subsequence
$k_j$ so rapidly that
$\norm{f_{k_{j+1}}-f_{k_j}}_1\leq2^{-3j}$. Then
\[
 \ga\bigl(\{|\Rg f_{k_{j+1}}-\Rg f_{k_j}|>2^{-j}\}\bigr)
       \leq2^{1-2j}.
\]
The right-hand side is summable. By the Borel--Cantelli lemma, almost
everywhere all sufficiently late increments have norm at most
$2^{-j}$. Hence the subsequence converges almost everywhere to a finite
vector-valued function, denoted $Tf$. The Cauchy-in-measure property
implies that the full sequence converges to the same limit in measure.
The limit is independent of the approximating sequence: applying
\eqref{eq:Cauchy-measure} to the difference of two approximations shows
that their output difference tends to zero in measure.

For any $t>0$, almost-everywhere convergence along the chosen subsequence
gives
\[
 \one_{\{|Tf|>t\}}\leq
       \liminf_{j\to\infty}\one_{\{|\Rg f_{k_j}|>t\}}.
\]
Indeed, at a point with $|Tf|>t$ there is a strict positive margin above
$t$, and the approximating norms eventually exceed $t$. Fatou's lemma
and the bound already proved for all levels on $L^2$ inputs
(that is, \eqref{eq:small-level} together with \eqref{eq:main-bound})
give, using $\norm{f_{k_j}}_1\to\norm{f}_1$,
\begin{equation}\label{eq:limit-bound}
 \ga(\{|Tf|>t\})
 \leq\liminf_j\ga(\{|\Rg f_{k_j}|>t\})
 \leq\frac{2}{t}\norm{f}_1.
\end{equation}

Linearity follows by approximating two functions simultaneously and using
continuity of addition and scalar multiplication for convergence in
measure. The operator agrees with \cref{def:riesz} on $L^2$, by using
the constant approximating sequence for an $L^2$ function. We denote it
again by $\Rg$. Applying \eqref{eq:limit-bound} to differences gives
\[
 \norm{\Rg f-\Rg h}_{L^{1,\infty}(\ga;\R^n)}
       \leq2\norm{f-h}_1,
\]
so the extension is continuous in the weak-$L^1$ quasi-norm. Uniqueness
among continuous linear extensions follows from the density of $L^2$
in $L^1$. This completes the proof of \cref{thm:endpoint}.

\begin{remark}[A mean-sensitive bound at large levels]
\label{rem:mean-tail}
Mass preservation in the two positive decompositions gives
$\me(g)=\me(f)$. Thus
$\norm{\Pi_0g}_2^2=\norm{g}_2^2-\me(f)^2$ in the proof above,
and retaining this term yields
\begin{equation}\label{eq:mean-tail}
 t\,\ga(\{|\Rg f|>t\})
 \leq2\norm{f}_1-\frac{\me(f)^2}{t},
 \qquad t>\norm{f}_1.
\end{equation}
This first holds on $L^2$. For a fixed $t>\norm{f}_1$ and an
$L^1$ approximation, eventually $t>\norm{f_k}_1$; convergence of the
means and the same Fatou argument prove it for every real $f\in L^1$.
\end{remark}

\subsection{Identification with the principal-value transform}
\label{subsec:principal-value}

The argument so far concerns the spectral transform. We now identify
its $L^1$ extension with the kernel realization. Equality of off-diagonal
kernels alone would not suffice, since it would not detect a term
supported on the diagonal. The cancellation on smooth inputs is checked
in \cref{app:pv-core}; the passage to $L^1$ uses only fixed-dimensional
maximal estimates.

\begin{proposition}[Principal-value realization]
\label{prop:pv-identification}
For every $f\in L^1(\ga)$ and $j=1,\ldots,n$, the principal value
\begin{equation}\label{eq:pv-realization}
 \lim_{\varepsilon\downarrow0}
 \int_{|x-y|>\varepsilon}K_{\gamma,j}(x,y)f(y)\,dy
\end{equation}
exists for almost every $x$ and equals the extension $R_{\gamma,j}f$
constructed in \cref{subsec:L1extension}. Here
$K_{\gamma,j}$ is the \emph{Lebesgue} kernel in
\eqref{eq:normalized-kernel}. Consequently, the same full-vector
bound $2$ holds for this principal-value realization on real inputs.
\end{proposition}

\begin{proof}
Fix $n$ and $j$, and denote the integral in
\eqref{eq:pv-realization} by $T_{j,\varepsilon}f$. The
fixed-dimensional theorem of Fabes, Guti\'errez, and Scotto
\cite{FGS94}, with the multiplicative normalization in
\eqref{eq:normalized-kernel}, gives
\begin{equation}\label{eq:fixed-n-maximal}
 \ga\left(\left\{\sup_{\varepsilon>0}|T_{j,\varepsilon}h|>t\right\}\right)
 \leq\frac{C_{n,j}}t\norm{h}_1,\qquad t>0.
\end{equation}
No uniformity in $n$ is needed in this step. By \cref{lem:pv-core}, the limit agrees with
$R_{\gamma,j}\phi$ for every $\phi\in C_c^\infty(\R^n)$.

Choose $f_k\in C_c^\infty$ with $f_k\to f$ in $L^1(\ga)$.
For the oscillation at zero,
\[
 \omega_f(x):=\lim_{\rho\downarrow0}
 \sup_{0<\varepsilon,\delta<\rho}
 |T_{j,\varepsilon}f(x)-T_{j,\delta}f(x)|,
\]
linearity and convergence for $f_k$ imply
\[
 \omega_f\leq2\sup_{\varepsilon>0}
                  |T_{j,\varepsilon}(f-f_k)|.
\]
Thus \eqref{eq:fixed-n-maximal} gives
$\ga(\{\omega_f>t\})\leq2C_{n,j}\norm{f-f_k}_1/t$.
Letting $k\to\infty$ shows that $\omega_f=0$ almost everywhere.
The maximal function in \eqref{eq:fixed-n-maximal} is finite almost
everywhere as well, so the principal value exists and is finite.
Call it $T_jf$. Passing to the limit in each truncation, or using the
pointwise bound by the maximal operator, gives
\[
 \ga(\{|T_jf-T_jf_k|>t\})
 \leq C_{n,j}\norm{f-f_k}_1/t.
\]
Consequently $T_jf_k\to T_jf$ in measure. The spectral extension also
satisfies $R_{\gamma,j}f_k\to R_{\gamma,j}f$ in measure. Since the
limits agree on every $f_k$, uniqueness of a limit in measure gives
$T_jf=R_{\gamma,j}f$ almost everywhere. There are finitely many
coordinates in the fixed dimension, so the vector realizations agree
almost everywhere. The constant $2$ is supplied by
\cref{thm:endpoint}, not by the constants $C_{n,j}$ used for this
identification.
\end{proof}

\subsection{Coordinate estimates, normalization, and basic checks}

\begin{corollary}[Coordinate and historical spectral conventions]
For every $j\in\{1,\ldots,n\}$ and real $f\in L^1(\ga)$,
\[
 \sup_{t>0}t\ga(\{|R_{\gamma,j}f|>t\})\leq2\norm{f}_1.
\]
For the vector operator $\mathcal R=-\nabla L^{-1/2}\Pi_0$ in
\eqref{eq:FGS-spectral},
\begin{equation}\label{eq:historical-bound}
 \sup_{t>0}t\ga(\{|\mathcal Rf|>t\})
        \leq2\sqrt2\norm{f}_1.
\end{equation}
Both statements remain unchanged if $\ga$ is replaced by the
unnormalized measure $e^{-|x|^2}\,dx$, with the same operator.
\end{corollary}

\begin{proof}
Pointwise, $|R_{\gamma,j}f|\leq|\Rg f|_{\R^n}$.
The second estimate follows from $\mathcal R=-\sqrt2\Rg$ and the
homogeneity of the weak-$L^1$ quasi-norm. Scaling the measure scales both
sides by the same factor.
\end{proof}

Let $C_\gamma(n)$ denote the optimal real-input constant in
\cref{thm:endpoint}. The transform annihilates constants: $\Rg1=0$. For
$f(x)=\sqrt2\,x_j$, we have $Lf=f$ and $\Rg f=e_j$. Since
$\norm{f}_1=\sqrt{2/\pi}$, the optimal constant satisfies
\begin{equation}\label{eq:lower-check}
 \sqrt{\frac\pi2}\leq C_\gamma(n)\leq2.
\end{equation}
This is a normalization check, not a claim that either bound is sharp.

The kernel prefactor also matters. Replacing the normalized kernels in
\eqref{eq:normalized-kernel} by $k_j^{\rm raw}$ rescales the transform to
\[
 T^{\rm raw}=\frac{\pi^{(n+1)/2}}{\sqrt2}\Rg.
\]
The same test then gives the lower bound
$\pi^{(n+2)/2}/2$ for its weak-type constant. Thus the dimension-free
estimate concerns the normalized spectral transform, not the raw kernel
family.

\subsection{Complex-valued inputs}
\label{subsec:complex}

\begin{corollary}[A complex-valued estimate]
\label{cor:complex}
The complex-linear extension satisfies
\begin{equation}\label{eq:complex-bound}
 \sup_{t>0}t\,\ga(\{|\Rg f|_{\mathbb C^n}>t\})
 \leq\frac8\pi\norm{f}_{L^1(\ga;\mathbb C)}.
\end{equation}
\end{corollary}

\begin{proof}
First let $f=f_1+if_2\in L^2$ and put
$S=\norm{f_1}_1+\norm{f_2}_1$. For $t>\norm{f}_1$, apply the
four nonnegative decompositions of $f_1^+,f_1^-,f_2^+,f_2^-$ at height
$t$. Every input has mean at most $\norm{f}_1<t$. If $g_1,g_2$ are
the signed good functions and $\Omega$ is the union of the four
positivity sets, then
\[
 \ga(\Omega)\leq S/t,\qquad
 \norm{g_1+ig_2}_2^2
 =\norm{g_1}_2^2+\norm{g_2}_2^2\leq tS.
\]
The complex $L^2$ identity follows by applying \eqref{eq:R-L2} to
real and imaginary parts. The transformed bad part vanishes outside
$\Omega$, so the same Chebyshev calculation yields
\[
 t\,\ga(\{|\Rg f|_{\mathbb C^n}>t\})\leq2S.
\]
Rather than use $S\leq\sqrt2\norm{f}_1$, repeat this argument for
$e^{i\theta}f$. Complex linearity implies
$|\Rg(e^{i\theta}f)|=|\Rg f|$. For each $\theta$, the right-hand
side becomes $2S_\theta$, where
\[
 S_\theta=\norm{\operatorname{Re}(e^{i\theta}f)}_1
          +\norm{\operatorname{Im}(e^{i\theta}f)}_1.
\]
Tonelli's theorem and a rotation of the scalar complex number $f(x)$
give
\[
 \frac1{2\pi}\int_0^{2\pi}S_\theta\,d\theta
 =\frac4\pi\norm{f}_1.
\]
Averaging the preceding bound proves \eqref{eq:complex-bound} for
$t>\norm{f}_1$. Smaller levels are bounded by $\norm{f}_1$ since
$\ga$ has mass one. Approximation by complex $L^2$ functions, with the
Cauchy-in-measure and Fatou arguments of \cref{subsec:L1extension},
then gives a unique continuous complex-linear extension to $L^1$.
\end{proof}

\subsection{Maximal truncations and the scope of the result}

For an off-diagonal vector kernel $K(x,y)$, its hard truncation at radius
$\varepsilon>0$ and its maximal truncation are defined by
\begin{equation}\label{eq:maximal-truncation}
 T_\varepsilon f(x)
   =\int_{|x-y|>\varepsilon}K(x,y)f(y)\,dy,
 \qquad
 T_*f(x)=\sup_{\varepsilon>0}|T_\varepsilon f(x)|_{\R^n}.
\end{equation}
For the Gaussian transform, $K=K_\gamma$ is the vector kernel whose
coordinates $K_{\gamma,j}$ are given in \eqref{eq:normalized-kernel};
for the standard Euclidean vector Riesz transform,
$K(x,y)=c_n(y-x)/|x-y|^{n+1}$ for the convention
$\nabla(-\Delta)^{-1/2}$, where
$c_n=\Gamma((n+1)/2)/\pi^{(n+1)/2}$.

The argument here proves an estimate for the untruncated operator
$\Rg$, both spectrally and in its principal-value realization. Its localization step uses the exact identity
$\Rg L^{1/2}u=Du$. Hard truncation does not preserve this identity, so the
obstacle argument does not yield a weak-type $(1,1)$ estimate for the
Gaussian maximal truncation $R_{\gamma,*}$ defined by
\eqref{eq:maximal-truncation}. The same limitation applies to the
maximal truncation of the standard Euclidean Riesz transform. In
particular, the dimension-free weak-type result for the untruncated
Euclidean transform in \cite{OSS26} does not settle its maximal-truncation
analogue. A dimension-independent weak-type bound for either maximal truncation
is a separate question, not a consequence of the present argument. The fixed-dimensional maximal-truncation estimates for the
Gaussian transforms in \cite{FGS94} do not provide such a
dimension-independent bound. For the Euclidean maximal transforms,
dimension-free strong $L^p$ estimates for $1<p<\infty$ are available
\cite{KW23,KWZ26}; these do not by themselves imply the weak endpoint.

\appendix
\section{Spectral and kernel agreement on smooth inputs}
\label{app:pv-core}

We give the cancellation argument used in
\cref{prop:pv-identification}. Throughout this appendix $n$ and $j$ are
fixed. Constants denoted by $C_n$ may depend on $n$; they serve only to
identify two realizations of the same operator and do not enter the
endpoint estimate.

\begin{lemma}[Agreement on a smooth core]
\label{lem:pv-core}
For every $f\in C_c^\infty(\R^n)$ and every $x\in\R^n$,
\begin{equation}\label{eq:pv-core}
 \lim_{\varepsilon\downarrow0}
 \int_{|x-y|>\varepsilon}K_{\gamma,j}(x,y)f(y)\,dy
 =R_{\gamma,j}f(x),
\end{equation}
where the right-hand side is the continuous semigroup representative
in \eqref{eq:R-gradient-semigroup}.
\end{lemma}

\begin{proof}
\textit{Step 1: Isolate the odd part of the heat-kernel derivative.}
Fix $x$ and $s>0$, and write
\[
 r=e^{-s},\qquad \sigma=1-r^2,\qquad d=(1-r)x,\qquad z=y-x,
 \qquad q_\sigma(z)=(\pi\sigma)^{-n/2}e^{-|z|^2/\sigma}.
\]
Equation \eqref{eq:Mehler-Lebesgue} becomes
$p_s(x,y)=q_\sigma(z+d)$. Hence
\begin{equation}\label{eq:pv-derivative}
 D_j^x p_s(x,y)=\sqrt2\,r\,V_{j,\sigma}(z+d),
 \qquad V_{j,\sigma}(z)=\frac{z_j}{\sigma}q_\sigma(z).
\end{equation}
The function $V_{j,\sigma}$ is odd in $z_j$, and scaling gives
\begin{equation}\label{eq:pv-translation}
 \norm{\nabla V_{j,\sigma}}_{L^1(dx)}\leq C_n/\sigma,
 \qquad
 \norm{V_{j,\sigma}(\cdot+d)-V_{j,\sigma}}_{L^1(dx)}
 \leq C_n|d|/\sigma.
\end{equation}
For the second inequality use the fundamental theorem of calculus
along the segment from $z$ to $z+d$, followed by Fubini.

\medskip
\noindent\textit{Step 2: Check the principal-value cancellation on constants.}
Set
\[
 I_{\varepsilon,s}(x)
 :=\int_{|x-y|>\varepsilon}D_j^xp_s(x,y)\,dy.
\]
The integral of $V_{j,\sigma}$ over the centered exterior ball is zero.
Thus \eqref{eq:pv-translation} implies
\begin{equation}\label{eq:pv-constant-majorant}
 |I_{\varepsilon,s}(x)|
 \leq C_n r\frac{|d|}{\sigma}
 =C_n|x|\frac{e^{-s}}{1+e^{-s}}.
\end{equation}
This bound, independent of $\varepsilon$, is integrable against
$ds/\sqrt{s}$ on $(0,\infty)$. For each fixed $s$,
$I_{\varepsilon,s}(x)\to\int D_j^xp_s(x,y)\,dy=0$ as
$\varepsilon\downarrow0$.

For each fixed $\varepsilon>0$, Fubini is legitimate in the kernel
integral over $|x-y|>\varepsilon$. Indeed, at small $s$ we have
$|d|<\varepsilon/2$, so the change of variables $w=z+d$ gives a
Gaussian tail with $|w|>\varepsilon/2$; its exponential decay in
$\varepsilon^2/\sigma$ dominates the powers of $\sigma^{-1}$.
For $s$ bounded away from zero use
$\int|D_j^xp_s(x,y)|\,dy\leq C_n e^{-s}/\sqrt{1-e^{-2s}}$.
These bounds prove absolute integrability of the double integral for
fixed $x,\varepsilon$. Therefore
\[
 \int_{|x-y|>\varepsilon}K_{\gamma,j}(x,y)\,dy
 =\frac1{\sqrt\pi}\int_0^\infty
       I_{\varepsilon,s}(x)\,\frac{ds}{\sqrt s}.
\]
Dominated convergence using \eqref{eq:pv-constant-majorant} yields
\begin{equation}\label{eq:pv-constant-cancellation}
 \lim_{\varepsilon\downarrow0}
 \int_{|x-y|>\varepsilon}K_{\gamma,j}(x,y)\,dy=0.
\end{equation}
In particular, no diagonal multiplication term is left undetermined.

\medskip
\noindent\textit{Step 3: Subtract the value of the input at the center.}
Let $L_f=\norm{\nabla f}_\infty$. Since
$|f(y)-f(x)|\leq L_f|z|$, \eqref{eq:pv-derivative} and the
substitution $w=z+d$ give
\begin{align*}
 &\int |D_j^xp_s(x,y)|\,|f(y)-f(x)|\,dy\\
 &\quad\leq\frac{\sqrt2 rL_f}{\sigma}
       \int |w_j|(|w|+|d|)q_\sigma(w)\,dw\\
 &\quad\leq C_nL_f r\left(1+\frac{|d|}{\sqrt\sigma}\right)
 \leq C_nL_f e^{-s}(1+|x|).
\end{align*}
Here the Gaussian moments scale as
$\int |w_j||w|q_\sigma\leq C_n\sigma$ and
$\int |w_j|q_\sigma\leq C_n\sqrt\sigma$, and
$(1-r)/\sqrt{1-r^2}\leq1$. The final bound is integrable against
$ds/\sqrt{s}$. We can consequently interchange the time and space
integrals to obtain
\begin{align}
 \int K_{\gamma,j}(x,y)(f(y)-f(x))\,dy
 &=\frac1{\sqrt\pi}\int_0^\infty
        D_jP_sf(x)\,\frac{ds}{\sqrt s}\notag\\
 &=R_{\gamma,j}f(x).
 \label{eq:pv-subtracted}
\end{align}
The left-hand side is absolutely convergent; the equality uses
$\int D_j^xp_s(x,y)\,dy=0$ and
\eqref{eq:R-gradient-semigroup}.

Finally write each truncated integral as
\[
 \int_{|x-y|>\varepsilon}K_{\gamma,j}(x,y)(f(y)-f(x))\,dy
 +f(x)\int_{|x-y|>\varepsilon}K_{\gamma,j}(x,y)\,dy.
\]
Absolute convergence in \eqref{eq:pv-subtracted} treats the first
term, and \eqref{eq:pv-constant-cancellation} treats the second.
This proves \eqref{eq:pv-core}.
\end{proof}

\section*{Acknowledgments}

The authors gratefully acknowledge OpenAI for generously providing
them with complimentary academic accounts and for making available
the Prism \LaTeX{} editor, which greatly facilitated the preparation
of this manuscript. Discussions with ChatGPT were valuable in
exploring ideas, developing and organizing mathematical arguments,
and refining the exposition. An AI-generated report obtained using
Claude was also consulted in revising the manuscript.
The authors take full responsibility for all mathematical statements,
proofs, references, and conclusions presented in the paper.

\bigskip
\noindent\begin{minipage}{\textwidth}
\small
\textsc{Cristian E. Guti\'errez}\\
Department of Mathematics, Temple University\\
Email: \href{mailto:gutierre@temple.edu}{\texttt{gutierre@temple.edu}}

\medskip
\textsc{Juan Manfredi}\\
Department of Mathematics, University of Pittsburgh\\
Email: \href{mailto:manfredi@pitt.edu}{\texttt{manfredi@pitt.edu}}
\end{minipage}

\end{document}